\documentclass[12pt,a4paper]{article}
\usepackage{amsfonts,amsmath,mathrsfs,amssymb, indentfirst}
\usepackage[dvips]{graphicx}
\usepackage{epsfig}
\usepackage{rotating}
\newtheorem{theorem}{Theorem}[section]
\newtheorem{lemma}[theorem]{Lemma}

\newtheorem{defin}{Definition}[section]
\newtheorem{examp}{Example}[section]
\newtheorem{remark}{Remark}[section]
\newtheorem{cor}[theorem]{Corollary}

\begin{document}

\title {Injective and non-injective realizations with symmetry}

\author{Bernd Schulze\\ Department of Mathematics and
Statistics, York University,\\ 4700 Keele Street, Toronto, ON M3J
1P3, Canada}

\maketitle

\begin{abstract}
In this paper, we introduce a natural classification of bar and joint frameworks that possess symmetry. This classification establishes the mathematical foundation for extending a variety of results in rigidity, as well as infinitesimal or static rigidity, to frameworks that are realized with certain symmetries and whose joints may or may not be embedded injectively in the space. In particular, we introduce a symmetry-adapted notion of `generic' frameworks with respect to this classification and show that `almost all' realizations in a given symmetry class are generic and all generic realizations in this class share the same infinitesimal rigidity properties.
Within this classification we also clarify under what conditions group representation theory techniques can be applied to further analyze the rigidity properties of a (not necessarily injective) symmetric realization.
\end{abstract}

\section{Introduction}

A $d$-dimensional bar and joint framework is a pair $(G,p)$, where $G$ is a graph and $p$ is a map that assigns to each vertex of $G$ a point in Euclidean $d$-space. A framework provides a mathematical model for a system where specified distances are maintained or a physical structure that consists of rigid bars
connected together at their ends by flexible joints.
Such structures, specifically the rigidity properties of such structures, are widely studied in a variety of sciences such as civil engineering (\cite{graver}), mechanical engineering (\cite{FG3}), and biochemistry (\cite{W5}). In each of these fields, symmetric structures are often of particular interest, since symmetry is not only exhibited by many objects found in nature (such as biomolecules, for example), but it is also a familiar feature of human built structures.\\\indent Some engineers and chemists have used group representation theory techniques to gain insight into the rigidity properties of symmetric structures \cite{FGsymmax, FG1, FG2, FG3, FG4, KG1, KG2, KG3}. While a variety of interesting and useful observations resulted from this approach, two kinds of shortcomings can be identified.\\\indent First, many of the results given in these references are not presented with a mathematically precise background or formulation nor with mathematically rigorous proofs.\\\indent The present paper therefore aims to provide the mathematical foundation that is necessary to give rigorous proofs of these results, as well as additional results and conjectures relating to the rigidity of symmetric frameworks (such as the ones stated in \cite{cfgsw, BS4, BS1, BS3, BS2}, for example).\\\indent This is accomplished by establishing the necessary mathematical terms and definitions and by introducing an appropriate classification of symmetric frameworks. One of the key properties of this classification is that for any given class, `almost all' realizations in this class share the same infinitesimal rigidity properties. This classification therefore not only sets the groundwork for symmetrizing results in rigidity, infinitesimal rigidity, and static rigidity, but it can also be used to develop a symmetry-adapted version of generic rigidity theory.\\\indent
The definitions and results presented in this paper have already been used to establish symmetrized versions of a variety of famous theorems in each of the above-mentioned theories (such as Maxwell's rule from 1864 and Laman's Theorem, for example) and are fundamental to the results in \cite{BS4, BS1, BS3, BS2}. \\\indent
Secondly, in their studies of symmetric frameworks, engineers and chemists have restricted their attention to frameworks whose joints are embedded injectively in the space. While this is a reasonable assumption for most applications (atoms of biomolecules or joints of 3-dimensional physical structures never coincide, for example), there are occasions where we do want to analyze frameworks with non-injective configurations (if we want to model a linkage in the plane with overlapping joints, for example). So, in order to obtain more general mathematical results and a more complete theory, we develop the mathematical foundation for the rigidity of symmetric structures in such a way that it also allows us to analyze symmetric frameworks with non-injective configurations.\\\indent
If one wants to apply group representation theory techniques to the analysis of non-injective symmetric realizations a variety of subtle difficulties can occur which we will address in Sections 5 and 6. We found a first indication of these difficulties  in probing the background for the symmetry-extended version of Maxwell's rule given in \cite{FGsymmax}. In examining some simple but extreme examples, we discovered that this rule (in its current version) does not give correct results when applied to some non-injective realizations. Therefore, in order to obtain a more general and still mathematically correct rule, it is necessary to reformulate this rule based on the mathematical foundation we establish in this paper.\\\indent
The structure of the paper is as follows:
in Section 2, we introduce mathematically explicit definitions of the relevant terms relating to symmetric structures that are frequently used in the chemistry (and engineering) literature. We also briefly explain the commonly used Schoenflies notation for point groups in dimensions 2 and 3, as we will be using this notation for all the examples throughout this paper.\\\indent In Section 3, our classification of symmetric frameworks is presented, along with some examples that will turn out to be very useful to illustrate some important observations in the following sections.\\\indent
Section 4 is devoted to establish a symmetry-adapted notion of a generic framework. We confirm that this new definition of generic satisfies all the desired properties and we also examine some additional important questions regarding this definition. The results in \cite{BS3} and many of the results in \cite{BS4} are based on this definition of generic.\\\indent In Sections 5 and 6, we illustrate the complications that may arise if one wants to extend rigidity results to non-injective symmetric realizations. More precisely, in Section 5, we examine how many distinct symmetry classes a given framework can possibly belong to, and in Section 6, we investigate under what conditions group representation theory techniques can be applied to the frameworks in a given symmetry class.\\\indent Finally, in Section 7, we briefly discuss possible extensions of our results to more general symmetric structures.


\section{Definitions and preliminaries}

We begin by introducing the necessary terms and definitions relating to symmetric frameworks. An introduction to infinitesimal rigidity is given in Section 4.

\begin{defin}\emph{A \emph{graph} $G$ is a finite
nonempty set of objects called \emph{vertices} together with a (possibly empty) set of unordered pairs of distinct vertices of
$G$ called \emph{edges}. The \emph{vertex set} of $G$ is denoted by $V(G)$ and the \emph{edge set} of $G$ is denoted by $E(G)$. Two vertices $u \ne v$ of a graph $G$ are said to be \emph{adjacent} if $\{u,v\}\in E(G)$, and \emph{independent} otherwise.}
\end{defin}

\begin{defin}\emph{An \emph{automorphism} of a graph $G$ is a permutation $\alpha$ of $V(G)$ such that $\{u,v\}\in E(G)$ if and only if $\{\alpha(u),\alpha(v)\}\in E(G)$.}
\end{defin}

The automorphisms of a graph $G$ form a group under composition which is denoted by $\textrm{Aut}(G)$.

\begin{defin}
\label{framework}
\emph{\cite{gss, W1, W2} A \emph{framework} (in $\mathbb{R}^{d}$) is a pair $(G,p)$, where $G$ is a graph and $p: V(G)\to \mathbb{R}^{d}$ is a map with the property that $p(u) \neq p(v)$ for all $\{u,v\} \in E(G)$. We also say that $(G,p)$
is a $d$-dimensional \emph{realization} of the \emph{underlying graph} $G$.}
\end{defin}

\begin{defin}
\emph{Let $(G,p)$ be a framework in $\mathbb{R}^{d}$. A \emph{joint} of $(G,p)$ is an ordered pair $\big(v,p(v)\big)$, where $v \in V(G)$. A \emph{bar} of $(G,p)$ is an unordered pair $\big\{\big(u,p(u)\big),\big(v,p(v)\big)\big\}$ of joints of $(G,p)$, where $\{u,v\} \in E(G)$. We define $\|p(u)-p(v)\|$ to be the \emph{length} of the bar $\big\{\big(u,p(u)\big),\big(v,p(v)\big)\big\}$, where $\|p(u)-p(v)\|$ is defined by the canonical inner product on $\mathbb{R}^{d}$.}
\end{defin}

Note that we allow a map $p$ of a framework $(G,p)$ to be non-injective, that is, two distinct joints $\big(u,p(u)\big)$ and $\big(v,p(v)\big)$ of $(G,p)$ may be located at the same point $p(u)=p(v)$ in $\mathbb{R}^{d}$, provided that $u$ and $v$ are independent vertices of $G$. However, if $\{u,v\} \in E(G)$, then $p(u)\neq p(v)$, and hence every bar $\big\{\big(u,p(u)\big),\big(v,p(v)\big)\big\}$ of $(G,p)$ has a strictly positive length.

We now establish the concept of a symmetric framework and give mathematically precise definitions of terms relating to symmetry which might have different meanings in different contexts. In the literature about symmetric structures it is common  to systematize the notion of symmetry by introducing the concept of a symmetry operation and its corresponding symmetry element \cite{bishop, cotton, Hall}. We begin with our definitions of these terms.

First, recall that an \emph{isometry} of $\mathbb{R}^{d}$ is a map $x:\mathbb{R}^{d}\to \mathbb{R}^{d}$ such that $\|x(a)-x(b)\|=\|a-b\|$ for all $a,b\in \mathbb{R}^{d}$.

\begin{defin}
\label{symop}
\emph{Let $(G,p)$ be a framework in $\mathbb{R}^{d}$.
A \emph{symmetry operation} of $(G,p)$ is an isometry $x$ of $\mathbb{R}^{d}$ such that for some $\alpha\in \textrm{Aut}(G)$, we have
$x\big(p(v)\big)=p\big(\alpha(v)\big)$ for all $v\in V(G)$.}
\end{defin}

A symmetry operation $x$ of a framework $(G,p)$ carries $(G,p)$ into a framework $(G,x \circ p)$ which is `geometrically
indistinguishable' from $(G,p)$. In other words, up to the labeling of the vertices of the underlying graph $G$, the frameworks $(G,p)$ and $(G,x \circ p)$ are the same.

\begin{defin}
\label{symel}
\emph{Let $x$ be a symmetry operation of a framework $(G,p)$ in $\mathbb{R}^{d}$.
The \emph{symmetry element} corresponding to $x$ is the affine subspace $F_{x}$ of $\mathbb{R}^{d}$ which
consists of all points in $\mathbb{R}^{d}$ that are fixed by $x$.}
\end{defin}

\begin{figure}[htp]
\begin{center}
\includegraphics[clip]{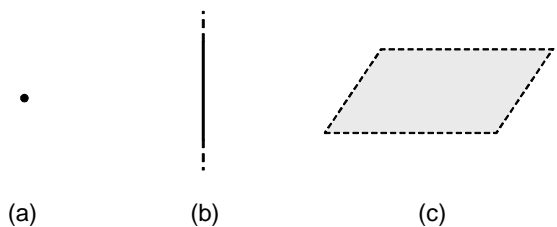}
\end{center}
\caption{\emph{Symmetry elements corresponding to symmetry operations in dimension $2$: (a) a rotation $C_{m}$, $m\geq 2$; (b) a reflection $s$; (c) the identity $Id$.}}
\label{symelem2}
\end{figure}

\begin{figure}[htp]
\begin{center}
\includegraphics[clip]{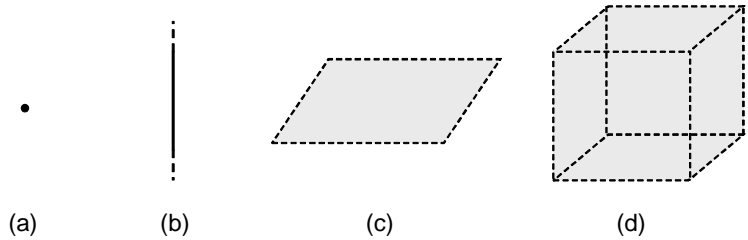}
\end{center}
\caption{\emph{Symmetry elements corresponding to symmetry operations in dimension $3$: (a) an improper rotation $S_{m}$, $m\geq 2$; (b) a rotation $C_{m}$, $m\geq 2$; (c) a reflection $s$; (d) the identity $Id$.}}
\label{symelem3}
\end{figure}

Since we only consider finite graphs, it follows directly from Definition \ref{symop} that a symmetry operation cannot be a translation. This implies in particular that a symmetry element is always non-empty. In fact, it is easy to see that if $x$ is a symmetry operation of a framework $(G,p)$ with $V(G)=\{v_{1},\ldots,v_{n}\}$, then the point $\frac{1}{n}\sum_{i=1}^{n}p(v_{i})$ must be fixed by $x$. Figures \ref{symelem2} and \ref{symelem3} depict all possible symmetry elements in dimensions 2 and 3.\\\indent
Note that distinct symmetry operations of a framework may have the same corresponding symmetry element. For example, distinct rotational symmetry operations of a 3-dimensional framework may share the same rotational axis.

The set of all symmetry operations of a given framework forms a group under composition. We adopt the following vocabulary from chemistry and crystallography:

\begin{defin} \emph{Let $(G,p)$ be a framework. Then the group which consists of all symmetry operations of $(G,p)$ is called the \emph{point group} of $(G,p)$.}
\end{defin}

For a systematic method to find the point group of a given framework, see \cite{bishop, cotton, Hall}, for example.

\begin{defin} \emph{A \emph{symmetry group} (in dimension $d$) is a subgroup of the orthogonal group $O(\mathbb{R}^{d})$ which consists of all orthogonal linear transformations of $\mathbb{R}^{d}$.}
\end{defin}

If $P$ is the point group of a $d$-dimensional framework, then, as noted above, there exists a point in
$\mathbb{R}^{d}$ which is fixed by all symmetry operations in $P$. Note that if the origin of $\mathbb{R}^{d}$ is fixed by a symmetry operation $x\in P$, then $x$ is an orthogonal linear transformation of $\mathbb{R}^{d}$. So, if the origin of $\mathbb{R}^{d}$ is fixed by every symmetry operation in $P$, then $P$ is a symmetry group.\\\indent Given a $d$-dimensional framework $(G,p)$, the framework $(G,T\circ p)$, where $T$ is a translation of $\mathbb{R}^d$, clearly has the same rigidity properties as $(G,p)$. Therefore, for our purposes we may wlog restrict our attention to frameworks whose point groups are symmetry groups. In this paper, the point group of every framework is assumed to be a symmetry group.

We use the Schoenflies notation to denote symmetry operations and symmetry groups in dimensions 2 and 3, as it is one of the standard notations in the literature about symmetric structures \cite{altherz, atkchil, bishop, cotton, Hall}.
\\\indent In the plane, the three kinds of possible symmetry operations are the identity $Id$, rotations $C_{m}$ about the origin by an angle of $\frac{2\pi}{m}$, where $m\geq 2$, and reflections $s$ in lines through the origin. The symmetry elements corresponding to these symmetry operations are shown in Figure \ref{symelem2}.

In the Schoenflies notation we differentiate between the following four types of symmetry groups in dimension 2: $\mathcal{C}_{1}$, $\mathcal{C}_{s}$, $\mathcal{C}_{m}$ and $\mathcal{C}_{mv}$, where $m\geq 2$.\\\indent $\mathcal{C}_{1}$ denotes the trivial group which only contains the identity $Id$. $\mathcal{C}_{s}$ denotes any symmetry group in dimension 2 that consists of the identity $Id$ and a single reflection $s$. For $m\geq 2$, $\mathcal{C}_{m}$ denotes any cyclic symmetry group of order $m$ which is generated by a rotation $C_{m}$, and $\mathcal{C}_{mv} $ denotes any symmetry group in dimension 2 that is generated by a pair $\{C_{m},s\}$.
\\\indent In 3-space, there are the following symmetry operations: the identity $Id$, rotations $C_{m}$ about axes through the origin by an angle of $\frac{2\pi}{m}$, where $m\geq 2$, reflections $s$ in planes through the origin, and improper rotations $S_{m}$ fixing the origin, where $m\geq 3$. An improper rotation $S_{m}$ is a rotation $C_{m}$ followed by the reflection $s$ whose symmetry element is the plane through the origin that is perpendicular to the axis of $C_{m}$. The axis of $C_{m}$ is called the improper rotation axis of $S_{m}$. By convention, $S_{1}$ and $S_{2}$ are treated separately, since $S_{1}$ is simply a reflection $s$ and $S_{2}$ is the inversion in the origin which is denoted by $i$.

This gives rise to the following families of possible symmetry groups in dimension 3: $\mathcal{C}_{1}$, $\mathcal{C}_{s}$, $\mathcal{C}_{i}$, $\mathcal{C}_{m}$, $\mathcal{C}_{mv}$, $\mathcal{C}_{mh}$, $\mathcal{D}_{m}$, $\mathcal{D}_{mh}$, $\mathcal{D}_{md}$, $\mathcal{S}_{2m}$, $\mathcal{T}$, $\mathcal{T}_{d}$, $\mathcal{T}_{h}$, $\mathcal{O}$, $\mathcal{O}_{h}$, $\mathcal{I}$, and $\mathcal{I}_{h}$, where $m\geq 2$. \\\indent
Analogous to the notation in dimension 2, $\mathcal{C}_{1}$ again denotes the trivial group that only contains the identity $Id$, $\mathcal{C}_{m}$ denotes any symmetry group in dimension 3 that is generated by a rotation $C_{m}$, where $m\geq 2$, and $\mathcal{C}_{s}$ denotes any symmetry group in dimension 3 that consists of the identity $Id$ and a single reflection $s$.\\\indent $\mathcal{C}_{i}$ is the symmetry group which consists of the identity and the inversion $i$ of $\mathbb{R}^{3}$.\\\indent $\mathcal{C}_{mv}$ denotes any symmetry group that is generated by a rotation $C_{m}$ and a reflection $s$ whose symmetry element contains the rotational axis of $C_{m}$. Similarly, a symmetry group $\mathcal{C}_{mh}$ is generated by a rotation $C_{m}$ and the reflection $s$ whose symmetry element is perpendicular to the axis of $C_{m}$.\\\indent The symbol $\mathcal{D}_{m}$ is used to denote a symmetry group in dimension 3 that is generated by a rotation $C_{m}$ and another 2-fold rotation $C_{2}$ whose rotational axis is perpendicular to the one of $C_{m}$. Symmetry groups of the types $\mathcal{D}_{mh}$ and $\mathcal{D}_{md}$ are generated by the generators $C_{m}$ and $C_{2}$ of a group $\mathcal{D}_{m}$ and by a reflection $s$. In the case of $\mathcal{D}_{mh}$, the symmetry element of $s$ is the plane that is perpendicular to the $C_{m}$ axis and contains the origin (and hence contains the rotational axis of $C_{2}$), whereas in the case of $\mathcal{D}_{md}$, the symmetry element of $s$ is a plane that contains the $C_{m}$ axis and forms an angle of $\frac{\pi}{m}$ with the $C_{2}$ axis (i.e., the symmetry element of $s$ bisects the angle between adjacent half-turn axes created by rotating the $C_{2}$ axis about the $C_{m}$ axis).\\\indent If a symmetry group $S$ in dimension 3 is generated by an improper rotation $S_{k}$, where $k$ is even, say $k=2m$, then $S$ is denoted by $\mathcal{S}_{2m}$.\\\indent The remaining seven types of symmetry groups in dimension 3 are related to the Platonic solids and are placed into three divisions: the tetrahedral groups $\mathcal{T}$, $\mathcal{T}_{d}$ and $\mathcal{T}_{h}$, the octahedral groups $\mathcal{O}$ and $\mathcal{O}_{h}$, and the icosahedral groups $\mathcal{I}$ and $\mathcal{I}_{h}$.\\\indent
$\mathcal{T}$, $\mathcal{O}$, and $\mathcal{I}$ denote, respectively, a symmetry group that consists of all \emph{rotational} symmetry operations of a regular tetrahedron, octahedron, and icosahedron, whereas $\mathcal{T}_{h}$, $\mathcal{O}_{h}$, and $\mathcal{I}_{h}$ denote, respectively, a symmetry group that consists of \emph{all} symmetry operations of a regular tetrahedron, octahedron, and icosahedron. $\mathcal{T}_{d}$ denotes a symmetry group that is generated by the elements of a group $\mathcal{T}$ and those three reflections whose symmetry elements each contain two of the three axes that correspond to half-turns in $\mathcal{T}$.\\\indent For a more detailed description of these groups we refer the reader to \cite{altherz, atkchil, bishop, cotton, Hall, BS4}, for example.


\section{A classification of symmetric frameworks}

In order to symmetrize results in rigidity theory, particularly results in generic rigidity theory, we need an appropriate classification of symmetric frameworks. Naturally, we require that frameworks in the same class have the same underlying graph. This classification should also be such that `almost all' frameworks within a given class share the same infinitesimal rigidity properties, so that we can develop a symmetrized version of generic rigidity theory with respect to this classification.

\begin{defin}
\label{symclass}
\emph{Let $G$ be a graph and $S$ be a symmetry group in dimension $d$. Then we denote
\emph{$\mathscr{R}_{(G,S)}$} to be the set of all $d$-dimensional realizations of $G$ whose point group is either equal to $S$ or contains $S$ as a subgroup. For an element of $\mathscr{R}_{(G,S)}$ we say that it is a \emph{realization of the pair $(G,S)$}.}
\end{defin}

\begin{theorem}
\label{symcharac}
Let $(G,p)$ be a $d$-dimensional realization of a graph $G$ and $S$ be a symmetry group in dimension $d$. Then $(G,p)\in \mathscr{R}_{(G,S)}$ if and only if there exists a map $\Phi:S\to \textrm{Aut}(G)$ such that
$x\big(p(v)\big)=p\big(\Phi(x)(v)\big)$ for all $v\in V(G)$ and all $x\in S$.
\end{theorem}
\textbf{Proof.} It follows immediately from the definitions that $(G,p)\in \mathscr{R}_{(G,S)}$ if and only if $S$ is a subgroup of the point group of $(G,p)$ if and only if every element of $S$ is a symmetry operation of $(G,p)$ if and only if for every $x\in S$, there exists an automorphism $\alpha_{x}$ of $G$ that satisfies $x\big(p(v)\big)=p\big(\alpha_{x}(v)\big)$ for all $v\in V(G)$. $\square$

\begin{remark}
\emph{Note that a set of the form $\mathscr{R}_{(G,S)}$ can possibly be empty. For example, there clearly exists no realization of $(K_{2},\mathcal{C}_{3})$, where $K_{2}$ is the complete graph on $2$ vertices and $\mathcal{C}_{3}$ is a symmetry group in dimension 2.}
\end{remark}

Theorem \ref{symcharac} gives rise to the following natural classification of the frameworks within a set $\mathscr{R}_{(G,S)}$.

\begin{defin}
\label{classif} \emph{Let $S$ be a symmetry group, $(G,p)$ be a framework in $\mathscr{R}_{(G,S)}$ and $\Phi$ be a map from $S$ to $\textrm{Aut}(G)$. Then $(G,p)$ is said
to be \emph{of type $\Phi$} if the following equations hold:
\begin{equation}\label{eq:p} x\big(p(v)\big)=p\big(\Phi(x)(v)\big) \textrm{ for all } v \in V(G) \textrm{ and all } x\in S \textrm{.}\nonumber\end{equation}
We denote $\mathscr{R}_{(G,S,\Phi)}$ to be the set of all realizations of $(G,S)$ which are of type $\Phi$.}
\end{defin}

Given a graph $G$ and a symmetry group $S$ in dimension $d$, different choices of types $\Phi:S\to \textrm{Aut}(G)$ frequently lead to very different geometric types of realizations of $(G,S)$. This is because a type $\Phi$ forces the joints and bars of a framework in $\mathscr{R}_{(G,S,\Phi)}$ to assume certain geometric positions in $\mathbb{R}^d$. We give a few examples for small symmetry groups in dimensions 2 and 3 to demonstrate this.

\begin{figure}[htp]
\begin{center}
\includegraphics[clip]{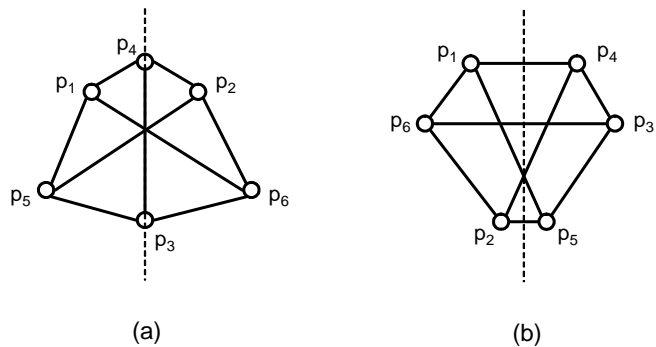}
\end{center}
\caption{\emph{$2$-dimensional realizations of $(K_{3,3},\mathcal{C}_{s})$ of different types.}}
\label{K33types}
\end{figure}

\begin{examp}
\label{K33ex}
\emph{Figure \ref{K33types} shows two realizations of $(K_{3,3},\mathcal{C}_{s})$ of different types, where $K_{3,3}$ is the complete bipartite graph with partite sets $\{v_{1},v_{2},v_{3}\}$ and $\{v_{4},v_{5},v_{6}\}$ and $\mathcal{C}_{s}=\{Id,s\}$ is a symmetry group in dimension 2 generated by a reflection. The framework in Figure \ref{K33types} (a) is a realization of $(K_{3,3},\mathcal{C}_{s})$ of type
$\Phi_{a}$, where $\Phi_{a}: \mathcal{C}_{s} \to \textrm{Aut}(K_{3,3})$ is defined by
\begin{eqnarray} \Phi_{a}(Id)& =& id\nonumber\\\Phi_{a}(s)&=&
(v_{1}\,v_{2})(v_{5}\,v_{6})(v_{3})(v_{4})\textrm{,}\nonumber
\end{eqnarray}
and the framework in Figure \ref{K33types} (b)
is a realization of $(K_{3,3},\mathcal{C}_{s})$ of type
$\Phi_{b}$, where $\Phi_{b}: \mathcal{C}_{s} \to \textrm{Aut}(K_{3,3})$ is defined by
\begin{eqnarray}\Phi_{b}(Id)& = &id\nonumber\\ \Phi_{b}(s) &=& (v_{1}\,v_{4})(v_{2}\,v_{5})(v_{3}\,v_{6})\textrm{.}\nonumber
\end{eqnarray}
Note that for any framework $(K_{3,3},p)$ in the set $\mathscr{R}_{(K_{3,3},\mathcal{C}_{s},\Phi_{a})}$, the points $p_{3}$ and $p_{4}$ must lie in the symmetry element corresponding to $s$ (i.e., in the mirror line of $s$), because $s\big(p(v_{i})\big)=p\big(\Phi_{a}(s)(v_{i})\big)=p(v_{i})$ for $i=3,4$. This says in particular that for any framework $(K_{3,3},p)$ in $\mathscr{R}_{(K_{3,3},\mathcal{C}_{s},\Phi_{a})}$, the entire undirected line segment $p_{3}p_{4}$ which corresponds to the bar $\big\{(v_{3},p_{3}),(v_{4},p_{4})\big\}$ of $(K_{3,3},p)$ must lie in the mirror line of $s$. We shall immediately become less formal and say that the bar $\big\{(v_{3},p_{3}),(v_{4},p_{4})\big\}$ lies in the mirror line of $s$.\\\indent Similarly, for any framework $(K_{3,3},p)$ in $\mathscr{R}_{(K_{3,3},\mathcal{C}_{s},\Phi_{b})}$, the bars $\big\{(v_{1},p_{1}),(v_{4},p_{4})\big\}$, $\big\{(v_{2},p_{2}),(v_{5},p_{5})\big\}$ and $\big\{(v_{3},p_{3}),(v_{6},p_{6})\big\}$ must be perpendicular to and centered at the mirror line of $s$.}
\end{examp}

\begin{figure}[htp]
\begin{center}
\includegraphics[clip]{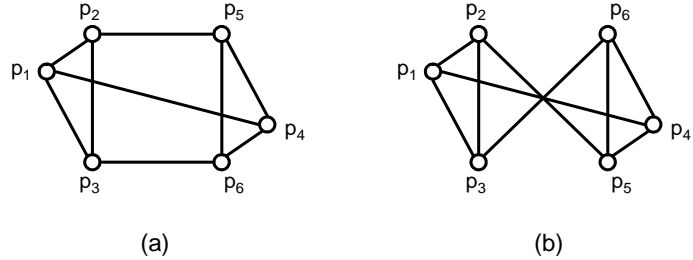}
\end{center}
\caption{\emph{$2$-dimensional realizations of $(G_{tp},\mathcal{C}_{2})$ of different types.}}
\label{triangprismtypes}
\end{figure}

\begin{examp}
\label{tpex}
\emph{Figure \ref{triangprismtypes} depicts two realizations of $(G_{tp},\mathcal{C}_{2})$ of different types, where $G_{tp}$ is the graph of a triangular prism and $\mathcal{C}_{2}=\{Id,C_{2}\}$ is the half-turn symmetry group in dimension 2. The framework in Figure \ref{triangprismtypes} (a) is a realization of $(G_{tp},\mathcal{C}_{2})$ of type $\Psi_{a}$, where $\Psi_{a}: \mathcal{C}_{2} \to \textrm{Aut}(G_{tp})$ is defined by
\begin{eqnarray} \Psi_{a}(Id)&=& id\nonumber\\
\Psi_{a}(C_{2}) & = &
(v_{1}\,v_{4})(v_{2}\,v_{6})(v_{3}\,v_{5})\textrm{.}\nonumber
\end{eqnarray}
and the framework in Figure \ref{triangprismtypes} (b) is a realization of $(G_{tp},\mathcal{C}_{2})$ of type $\Psi_{b}$, where $\Psi_{b}: \mathcal{C}_{2} \to \textrm{Aut}(G_{tp})$ is defined by
\begin{eqnarray} \Psi_{b}(Id) & = & id\nonumber\\ \Psi_{b}(C_{2}) & = &
(v_{1}\,v_{4})(v_{2}\,v_{5})(v_{3}\,v_{6})\textrm{.} \nonumber
\end{eqnarray}
It follows from the definitions of $\Psi_{a}$ and $\Psi_{b}$ that for any framework $(G_{tp},p)$ in $\mathscr{R}_{(G_{tp},\mathcal{C}_{2},\Psi_{a})}$, the bar $\big\{(v_{1},p_{1}),(v_{4},p_{4})\big\}$ must be centered at the origin (which is the center of the half-turn $C_{2}$), whereas for any framework $(G_{tp},p)$ in $\mathscr{R}_{(G_{tp},\mathcal{C}_{2},\Psi_{b})}$, all three bars $\big\{(v_{1},p_{1}),(v_{4},p_{4})\big\}$, $\big\{(v_{2},p_{2}),(v_{5},p_{5})\big\}$, and $\big\{(v_{3},p_{3}),(v_{6},p_{6})\big\}$ must be centered at the origin.}
\end{examp}

\begin{figure}[htp]
\begin{center}
\includegraphics[clip]{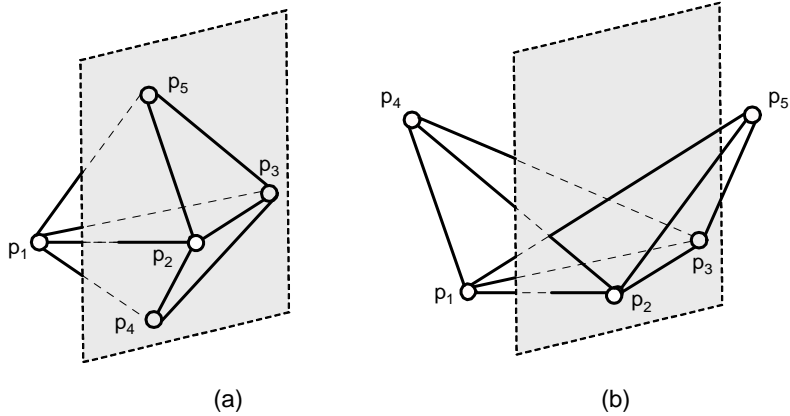}
\end{center}
\caption{\emph{$3$-dimensional realizations of $(G_{bp},\mathcal{C}_{s})$ of different types.}}
\label{bipyr}
\end{figure}

\begin{examp}
\label{3Dtypes}
\emph{Finally, Figure \ref{bipyr} depicts two realizations of $(G_{bp},\mathcal{C}_{s})$ of different types, where $G_{bp}$ is the graph of a triangular bipyramid and $\mathcal{C}_{s}=\{Id,s\}$ is a symmetry group in dimension 3. The framework in Figure \ref{bipyr} (a) is an element of $\mathscr{R}_{(G_{bp},\mathcal{C}_{s},\Xi_{a})}$, where
$\Xi_{a}: \mathcal{C}_{s} \to \textrm{Aut}(G_{bp})$ is defined by
\begin{eqnarray} \Xi_{a}(Id) & = & id \nonumber\\
 \Xi_{a}(s) & = & (v_{1}\,v_{2})(v_{3})(v_{4})(v_{5})\textrm{, }\nonumber
\end{eqnarray} and the framework in Figure \ref{bipyr} (b) is an element of
$\mathscr{R}_{(G_{bp},\mathcal{C}_{s},\Xi_{b})}$, where $\Xi_{b}: \mathcal{C}_{s} \to \textrm{Aut}(G_{bp})$ is defined by
\begin{eqnarray}\Xi_{b}(Id) & = & id\nonumber\\ \Xi_{b}(s) & = & (v_{1}\,v_{2})(v_{4}\,v_{5})(v_{3})\textrm{.}\nonumber
\end{eqnarray}
For any framework $(G_{tp},p)$ in $\mathscr{R}_{(G_{bp},\mathcal{C}_{s},\Xi_{a})}$ or $\mathscr{R}_{(G_{bp},\mathcal{C}_{s},\Xi_{b})}$, the bar $\big\{(v_{1},p_{1}),(v_{2},p_{2})\big\}$ must be perpendicular to and centered at the mirror plane of $s$. Further, for any framework $(G_{tp},p)$ in $\mathscr{R}_{(G_{bp},\mathcal{C}_{s},\Xi_{a})}$, the joints $(v_{i},p_{i})$, $i=3,4,5$, must lie in the mirror plane of $s$, whereas for a framework $(G_{tp},p)$ in $\mathscr{R}_{(G_{bp},\mathcal{C}_{s},\Xi_{b})}$, only the joint $(v_{3},p_{3})$ must have this property and the joints $(v_{4},p_{4})$ and $(v_{5},p_{5})$ must be mirror images of each other with respect to $s$.}
\end{examp}

\begin{remark}\emph{Given a non-empty set $\mathscr{R}_{(G,S)}$, it is possible that $\mathscr{R}_{(G,S,\Phi)}=\emptyset$ for some map $\Phi:S\to \textrm{Aut}(G)$.\\\indent Consider, for example, the non-empty set $\mathscr{R}_{(K_{2},\mathcal{C}_{2})}$, where $\mathcal{C}_{2}=\{Id,C_{2}\}$ is the half-turn symmetry group in dimension 2, and let $I:\mathcal{C}_{2}\to \textrm{Aut}(K_{2})$ be the map which sends both $Id$ and $C_{2}$ to the identity automorphism of $K_{2}$. If $(K_{2},p)\in\mathscr{R}_{(K_{2},\mathcal{C}_{2},I)}$, then both joints of $(K_{2},p)$ must be located at the origin (which is the center of $C_{2}$). This contradicts Definition \ref{framework} of a framework, and hence we have $\mathscr{R}_{(K_{2},\mathcal{C}_{2},I)}=\emptyset$.}
\end{remark}

We will see in the next section that `almost all' frameworks within a set of the form $\mathscr{R}_{(G,S,\Phi)}$ share the same infinitesimal rigidity properties.

\section{The notion of $(S,\Phi)$-generic}

\subsection{Introduction to infinitesimal rigidity}

We briefly recall the necessary terms and definitions relating to infinitesimal rigidity.

\begin{defin}
\label{infinmotion}
\emph{Let $(G,p)$ be a framework in $\mathbb{R}^{d}$ with $V(G)=\{v_{1},v_{2},\ldots, v_{n}\}$. An \emph{infinitesimal motion} of $(G,p)$ is a function $u: V(G)\to \mathbb{R}^{d}$ such that
\begin{equation}
\label{infinmotioneq}
\big(p(v_{i})-p(v_{j})\big)\cdot \big(u(v_{i})-u(v_{j})\big)=0 \quad\textrm{ for all } \{v_{i},v_{j}\} \in E(G)\textrm{.}\end{equation}
An infinitesimal motion $u$ of $(G,p)$ is \emph{trivial} if there exists a family of differentiable functions $P_{i}:[0,1]\to \mathbb{R}^{d}, \, i=1,2,\ldots, n$, with $P_{i}(0)=p(v_{i})$ for all $i$ and $\|P_{i}(t)-P_{j}(t)\|=\|p(v_{i})-p(v_{j})\|$ for all $t\in [0,1]$ and all $1\leq i < j\leq n$ such that $u(v_{i})=P'_{i}(0)$ for all $i$.\\\indent
$(G,p)$ is said to be \emph{infinitesimally rigid} if every infinitesimal motion of $(G,p)$ is trivial. Otherwise $(G,p)$ is said to be \emph{infinitesimally flexible}. See \cite{gss, W1}, for example, for more details.}
\end{defin}

Definition \ref{infinmotion} is motivated by the idea of a motion that displaces the joints of $(G,p)$ on differentiable displacement paths while preserving the length $\| p(v_{i})-p(v_{j})\|$ of all bars $\big\{\big(v_{i},p(v_{i})\big),\big(v_{j},p(v_{j})\big)\big\}$ of $(G,p)$. A simple differentiation yields the equations for the velocity vectors given in (\ref{infinmotioneq}) \cite{graver, gss}. So an infinitesimal motion of a framework $(G,p)$ is a set of initial velocity vectors, one at each joint, that neither stretch nor compress the bars of $(G,p)$. More precisely, condition (\ref{infinmotioneq}) says that for every edge $\{v_{i},v_{j}\} \in E(G)$, the projections of $u(v_{i})$ and $u(v_{j})$ onto the line through $p(v_{i})$ and $p(v_{j})$ have the same direction and the same length \cite{W1, W2}.

\begin{figure}[htp]
\begin{center}
\includegraphics[clip]{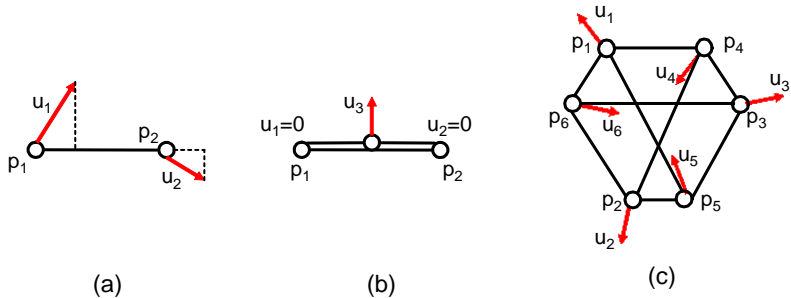}
\end{center}
\caption{\emph{The arrows indicate the non-zero velocity vectors of trivial (a) and non-trivial (b, c) infinitesimal motions of frameworks in $\mathbb{R}^2$.}}
\label{inmo}
\end{figure}

For a framework $(G,p)$ whose underlying graph $G$ has a vertex set that is indexed from 1 to $n$, say $V(G)=\{v_{1},v_{2},\ldots ,v_{n}\}$, we will frequently denote $p(v_{i})$ by $p_{i}$ for $i=1,2,\ldots, n$. The $k^{th}$ component of a vector $x$ we will denote by $(x)_{k}$.

The equations stated in Definition \ref{infinmotion} form a system of linear equations whose corresponding matrix is the so-called rigidity matrix. This matrix is fundamental in the study of infinitesimal rigidity \cite{graver, gss, W1, W2}.

\begin{defin}
\emph{ Let $G$ be a graph with $V(G)=\{v_{1},v_{2},\ldots,v_{n}\}$ and let $p:V(G)\to \mathbb{R}^{d}$. The
\emph{rigidity matrix} of $(G,p)$ is the $|E(G)| \times dn$ matrix  \\\indent \begin{displaymath} \mathbf{R}(G,p)=\left(
\begin{array} {ccccccccccc }
& & & & & \vdots & & & & & \\
0 & \ldots & 0 & p_{i}-p_{j}&0 &\ldots &0 & p_{j}-p_{i} &0 &\ldots &
0\\ & & & & & \vdots & & & & &\end{array}
\right)\textrm{,}\end{displaymath} i.e., for each edge $\{v_{i},v_{j}\}\in E(G)$, $\mathbf{R}(G,p)$ has the row with
$(p_{i})_{1}-(p_{j})_{1},\ldots,(p_{i})_{d}-(p_{j})_{d}$ in the columns $d(i-1)+1,\ldots,di$, $(p_{j})_{1}-(p_{i})_{1},\ldots,(p_{j})_{d}-(p_{i})_{d}$ in
the columns $d(j-1)+1,\ldots,dj$, and $0$ elsewhere.}
\end{defin}

\begin{remark}
\emph{Note that the rigidity matrix is defined for arbitrary pairs $(G,p)$, where $G$ is a graph and $p:V(G)\to \mathbb{R}^{d}$ is a map. If $(G,p)$ is not a framework, then there exists a pair of adjacent vertices of $G$ that are mapped to the same point in $\mathbb{R}^{d}$ under $p$ and every such edge of $G$ gives rise to a zero-row in $\mathbf{R}(G,p)$.}
\end{remark}

\begin{theorem}
\label{infinrigaff}\cite{asiroth}
A framework $(G,p)$ in $\mathbb{R}^d$ is infinitesimally rigid if and only if either $\textrm{rank }\big(\mathbf{R}(G,p)\big)=d |V(G)| - \binom{d+1}{2}$ or $G$ is a complete graph and the points $p(v)$, $v\in V(G)$, are affinely independent.
\end{theorem}

\begin{remark}
\label{affine}
\emph{Let $1\leq m\leq d$ and let $(G,p)$ be a framework in $\mathbb{R}^d$. If $(G,p)$ has at least $m+1$ joints and the points $p(v)$, $v\in V(G)$, span an affine subspace of $\mathbb{R}^{d}$ of dimension less than $m$, then $(G,p)$ is infinitesimally flexible. In particular, if $(G,p)$ is infinitesimally rigid and $|V(G)|\geq d$, then the points $p(v)$, $v\in V(G)$, span an affine subspace of $\mathbb{R}^d$ of dimension at least $d-1$. }
\end{remark}

\begin{defin}
\emph{Let $(G,p)$ be a framework. If the rows of the rigidity matrix $\mathbf{R}(G,p)$ are linearly independent, then $(G,p)$ is said to be \emph{independent}. The framework $(G,p)$ is said to be \emph{isostatic} if it is infinitesimally rigid and independent.}
\end{defin}

An isostatic framework has the property that it is minimal infinitesimally rigid, that is, it is infinitesimally rigid and the removal of any bar results in a framework that is not infinitesimally rigid \cite{gss, W1, W2}.

The main goal in combinatorial (or generic) rigidity is to establish infinitesimal rigidity properties that hold for `almost all' realizations $(G,p)$ of a graph $G$. The following standard definition of `generic' (\cite{graver, gss}) specifies what we mean by `almost all'.

\begin{defin}
\label{indetrigmatrix}
\emph{Let $K_{n}$ be the complete graph on $n$ vertices with $V(K_{n})=\{v_{1},v_{2},\ldots,v_{n}\}$. For each $i=1,2,\ldots ,n$, we introduce a $d$-tuple $p'_{i}=\big((p'_{i})_{1},\ldots, (p'_{i})_{d}\big)$ of variables and let
\begin{displaymath} \mathbf{R}(n,d) =\left(
\begin{array} {ccccccccccc }
& & & & & \vdots & & & & & \\
0 & \ldots & 0 & p'_{i}-p'_{j}&0 &\ldots &0 & p'_{j}-p'_{i} &0
&\ldots & 0\\ & & & & & \vdots & & & & &\end{array}
\right)\end{displaymath}
be the matrix that is obtained from the rigidity matrix $\mathbf{R}(K_{n},p)$ of a $d$-dimensional realization $(K_{n},p)$ by replacing each $(p_{i})_{j}\in\mathbb{R}$ with the variable $(p'_{i})_{j}$. We call $\mathbf{R}(n,d)$ the \emph{$d$-dimensional indeterminate rigidity matrix of} $K_{n}$.}
\end{defin}

\begin{defin}
\label{generic}
\emph{Let $V=\{v_{1},v_{2},\ldots,v_{n}\}$ and $p:V \to \mathbb{R}^d$ be a map. Further, let $K_{n}$ be the complete graph with $V(K_{n})=V$.\\\indent We say that $p$ is \emph{generic} if the determinant of any submatrix of $\mathbf{R}(K_{n},p)$ is zero only if the determinant of the corresponding submatrix of $\mathbf{R}(n,d)$ is (identically) zero.\\\indent
A framework $(G,p)$ is said to be \emph{generic} if $p$ is a generic map.}
\end{defin}

There are two fundamental facts regarding this definition of generic. First, the set of all non-generic maps $p$ of a finite set $V=\{v_{1},v_{2},\ldots,v_{n}\}$ to $\mathbb{R}^d$ is a closed set of measure zero. To see this, identify $p$ with a vector in $\mathbb{R}^{dn}$ and observe that the determinant of every submatrix of $\mathbf{R}(K_{n},p)$ is a polynomial in the variables $(p'_{i})_{j}$. If such a polynomial is not identically zero, then, by general algebraic geometry, it is non-zero for an open dense set of $p\in\mathbb{R}^{dn}$. Since $\mathbf{R}(K_{n},p)$ has only finitely many minors, the set of generic $p\in\mathbb{R}^{dn}$ is still an open dense subset of $\mathbb{R}^{dn}$.\\\indent Secondly, the infinitesimal rigidity properties are the same for all generic realizations of a graph $G$. This is specified in

\begin{theorem}
\label{genericrigtheorem}\cite{graver}
For a graph $G$ and a fixed dimension $d$, the following are equivalent:
\begin{itemize}
\item[(i)] $(G,p)$ is infinitesimally rigid (independent, isostatic) for some map $p:V(G)\to \mathbb{R}^d$;
\item[(ii)] every $d$-dimensional generic realization of $G$ is infinitesimally rigid (independent, isostatic).
\end{itemize}
\end{theorem}
This gives rise to

\begin{defin}
\emph{A graph $G$ is called \emph{generically infinitesimally rigid in dimension $d$} if $d$-dimensional generic realizations of $G$ are infinitesimally rigid. $G$ is called \emph{generically independent (isostatic) in dimension $d$} or \emph{generically $d$-independent ($d$-isostatic)} if $d$-dimensional generic realizations of $G$ are independent (isostatic).}
\end{defin}

\begin{remark}\emph{An easy but often useful observation concerning generic frameworks is that if a framework $(G,p)$ in $\mathbb{R}^d$ is generic, then the joints of $(G,p)$ are in \emph{general position}, that is, for $1\leq m\leq d$, no $m+1$ joints of $(G,p)$ lie in an $m-1$-dimensional affine subspace of $\mathbb{R}^d$.}
\end{remark}

\subsection{Infinitesimal rigidity for symmetric frameworks}

\begin{figure}[htp]
\begin{center}
\includegraphics[clip]{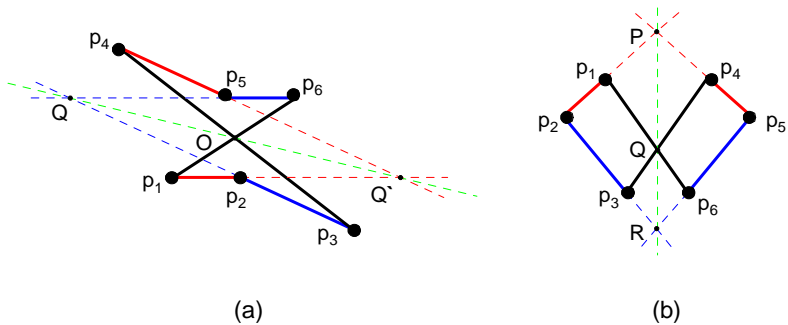}
\end{center}
\caption{\emph{By the converse of Pascal's Theorem, the joints of all realizations in $\mathscr{R}_{(K_{3,3},\mathcal{C}_{2})}$ and $\mathscr{R}_{(K_{3,3},\mathcal{C}_{s},\Phi_{b})}$ lie on a conic section.}}
\label{pas}
\end{figure}

Given a graph $G$, a non-trivial symmetry group $S$ in dimension $d$ clearly imposes restrictions on the possible geometric positions of realizations of $(G,S)$ in $\mathbb{R}^d$. In many cases, these restrictions can even be so strong that all realizations in a set of the form $\mathscr{R}_{(G,S)}$ are forced to be non-generic.\\\indent For example, every realization of $(K_{3},\mathcal{C}_{2})$, where $\mathcal{C}_{2}$ is a half-turn symmetry group in dimension 2 or 3, must be a degenerate triangle and is therefore non-generic. For a less trivial example, consider the complete bipartite graph $K_{3,3}$ and the symmetry group $\mathcal{C}_{2}$ in dimension 2. As shown in Figure \ref{pas} (a), the joints of any realization $(K_{3,3},p)$ in $\mathscr{R}_{(K_{3,3},\mathcal{C}_{2})}$ can be labeled in such a way that for the resulting hexagon $p_{1}\,p_{2}\,\ldots\,p_{6}$, there exists a pair of opposite sides (namely $p_{1}\,p_{6}$, $p_{3}\,p_{4}$) which intersect in the origin. If all three pairs of opposite sides of this hexagon are extended to their points of intersection, then the half-turn symmetry of $(K_{3,3},p)$ guarantees that these three points are collinear. Therefore, by the converse of Pascal's Theorem, the joints of $(K_{3,3},p)$ must lie on a conic section. It is well known that 2-dimensional realizations of $K_{3,3}$ whose joints lie on a conic section are in fact non-generic \cite{WW, W3}.\\\indent This shows that our notion of generic (without symmetry) is clearly not suitable any more once we restrict our attention to symmetric frameworks that lie within a set of the form $\mathscr{R}_{(G,S)}$.
\\\indent Note also that for a graph $G$, a symmetry group $S$, and two distinct maps $\Phi$ and $\Psi$ from $S$ to $\textrm{Aut}(G)$, it is possible that all realizations in $\mathscr{R}_{(G,S,\Phi)}$ are infinitesimally flexible, whereas `almost all' realizations in $\mathscr{R}_{(G,S,\Psi)}$ are isostatic.\\\indent For example, consider again the complete bipartite graph $K_{3,3}$, a symmetry group $\mathcal{C}_{s}$ in dimension 2, and the types $\Phi_{a}$ and $\Phi_{b}$ from Example \ref{K33ex}. $K_{3,3}$ is known to be a generically 2-isostatic graph and the pure condition (see \cite{WW}) for $K_{3,3}$ says that a 2-dimensional realization $(K_{3,3},p)$ is infinitesimally flexible if and only if the joints of $(K_{3,3},p)$ lie on a conic section. It follows (again from the converse of Pascal's Theorem) that every realization in $\mathscr{R}_{(K_{3,3},\mathcal{C}_{s},\Phi_{b})}$ is infinitesimally flexible (see also Figure \ref{pas} (b)), whereas `almost all' realizations in $\mathscr{R}_{(K_{3,3},\mathcal{C}_{s},\Phi_{a})}$ are isostatic.

Therefore, in order to define a modified, symmetry-adapted notion of generic for a set $\mathscr{C}\subseteq \mathscr{R}_{(G,S)}$ of symmetric frameworks in such a way that `almost all' realizations within $\mathscr{C}$ are generic and all generic realizations within $\mathscr{C}$ share the same infinitesimal rigidity properties, we need to restrict $\mathscr{C}$ to a set of the form $\mathscr{R}_{(G,S,\Phi)}$.

Let $G$ be a graph with $V(G)=\{v_{1},v_{2},\ldots,v_{n}\}$, $S$ be a symmetry group in dimension $d$, and $\Phi$ be a map from $S$ to $\textrm{Aut}(G)$. We will define a symmetry-adapted notion of generic for the set $\mathscr{R}_{(G,S,\Phi)}$ in an analogous way as we defined generic in Definition \ref{generic}. This requires the definition of a symmetry-adapted indeterminate rigidity matrix for $\mathscr{R}_{(G,S,\Phi)}$. The following observations set the groundwork for the definition of such a matrix.

Recall that for every framework $(G,p)$ in the set $\mathscr{R}_{(G,S,\Phi)}$, the equations stated in Definition \ref{classif} are satisfied, that is, we have $x\big(p(v_{i})\big)=p\big(\Phi(x)(v_{i})\big)$ for all $i=1,2,\ldots, n$ and all $x\in S$. Since every element of $S$ is an orthogonal linear transformation, we may identify each $x\in S$ with its corresponding orthogonal matrix $M_{x}$ that represents $x$ with respect to the canonical basis of $\mathbb{R}^{d}$. Therefore, for each $x\in S$, the equations in Definition \ref{classif} corresponding to $x$ form a system of linear equations which can be written as
\begin{displaymath} \mathbf{M^{(x)}} \left( \begin{array} {c} p_{1} \\ p_{2}
\\ \vdots \\ p_{n} \end{array} \right) =  \mathbf{P_{\Phi(x)}} \left( \begin{array} {c} p_{1} \\ p_{2}
\\ \vdots \\ p_{n} \end{array} \right) \textrm{, where}\end{displaymath}
\begin{displaymath} \mathbf{M^{(x)}}=\left( \begin{array} {cccc}
M_{x} & 0 & \ldots & 0 \\ 0 & M_{x} & \ddots & \vdots\\
\vdots & \ddots & \ddots & 0\\ 0 & \ldots &
0& M_{x} \end{array} \right)\textrm{,}\end{displaymath}
and $\mathbf{P_{\Phi(x)}}$ is the $dn \times dn$ matrix which is
obtained from the permutation matrix corresponding to $\Phi(x)$
by replacing each $1$ by a $d \times d$ identity matrix and each
$0$ by a $d \times d$ zero matrix. Equivalently, we have
\begin{displaymath} \Big(\mathbf{M^{(x)}}- \mathbf{P_{\Phi(x)}}\Big) \left( \begin{array} {c}
p_{1} \\ p_{2} \\ \vdots \\ p_{n} \end{array} \right) =
\mathbf{0}\textrm{.}
\end{displaymath}
We denote $L_{x,\Phi}=\textrm{ker }\big(\mathbf{M^{(x)}}- \mathbf{P_{\Phi(x)}}\big)$ and $U= \bigcap_{x \in S} L_{x,\Phi}$. Then $U$ is a subspace of $\mathbb{R}^{dn}$ which may be interpreted as the space of all those (possibly non-injective) configurations of $n$ points in $\mathbb{R}^{d}$ that possess the symmetry imposed by $S$ and $\Phi$. In particular, if we identify $p$ with a vector in $\mathbb{R}^{dn}$ (by using the order on $V(G)$), then $p$ is an element of $U$ whenever $(G,p)\in\mathscr{R}_{(G,S,\Phi)}$. Therefore, if we fix a basis $\mathscr{B}_{U}=\{u_{1},u_{2},\ldots,u_{k}\}$ of $U$, then every framework $(G,p) \in
\mathscr{R}_{(G,S,\Phi)}$ can be represented uniquely by the $k\times 1$ coordinate vector of $p$ relative to $\mathscr{B}_{U}$.

We are now ready to define the symmetry-adapted indeterminate rigidity matrix for $\mathscr{R}_{(G,S,\Phi)}$.

\begin{defin}
\label{indetsymrigmatrix} \emph{Let $G$ be a graph with $V(G)=\{v_{1},v_{2},\ldots,v_{n}\}$, $K_{n}$ be the complete graph with $V(K_{n})=V(G)$, $S$ be a symmetry group in dimension $d$, and $\Phi$ be a map from $S$ to $\textrm{Aut}(G)$. Further, let $\mathscr{B}_{U}=\{u_{1},u_{2},\ldots,u_{k}\}$ be a basis of $U=\bigcap_{x \in S} L_{x,\Phi}$. The \emph{symmetry-adapted indeterminate rigidity matrix} for $\mathscr{R}_{(G,S,\Phi)}$ (corresponding to $\mathscr{B}_{U}$) is the matrix
$\mathbf{R}_{\mathscr{B}_{U}}(n,d)$ which is obtained from the indeterminate rigidity matrix $\mathbf{R}(n,d)$ by introducing a $k$-tuple $(t'_{1},t'_{2},\ldots , t'_{k})$ of variables and replacing the $dn$ variables $(p'_{i})_{j}$ of $\mathbf{R}(n,d)$ as follows.\\ \indent
For each $i=1,2, \ldots ,n$ and each $j=1,\ldots,d$, we replace the variable $(p'_{i})_{j}$ in $\mathbf{R}(n,d)$ by the linear combination
$t'_{1}(u_{1})_{i_{j}}+t'_{2}(u_{2})_{i_{j}}+ \ldots + t'_{k}(u_{k})_{i_{j}}$.}
\end{defin}

\begin{remark}
\label{symindetrigmatrem}
\emph{Let $(G,p)\in \mathscr{R}_{(G,S,\Phi)}$ and $\mathscr{B}_{U}=\{u_{1},u_{2},\ldots,u_{k}\}$ be a basis of
$U=\bigcap_{x \in S} L_{x,\Phi}$. Then \begin{displaymath}\left( \begin{array} {c} p_{1} \\ p_{2} \\ \vdots \\
p_{n}\end{array} \right)=t_{1}u_{1}+\ldots+ t_{k}u_{k}\textrm{, for some }t_{1},\ldots,t_{k}\in \mathbb{R}\textrm{.}\end{displaymath}
So, if for $i=1,\ldots,k$, the variable $t'_{i}$ in $\mathbf{R}_{\mathscr{B}_{U}}(n,d)$ is replaced by
$t_{i}$ then we obtain the rigidity matrix $\mathbf{R}(K_{n},p)$ of the framework $(K_{n},p)$.}
\end{remark}

With the help of Definition \ref{indetsymrigmatrix} we can now also give the formal definition of our symmetry-adapted notion of generic for a set $\mathscr{R}_{(G,S,\Phi)}$.

\begin{defin}
\label{symgeneric}
\emph{Let $G$ be a graph with $V(G)=\{v_{1},v_{2},\ldots,v_{n}\}$, $K_{n}$ be the complete graph with $V(K_{n})=V(G)$, $S$ be a symmetry group in dimension $d$, $\Phi$ be a map from $S$ to $\textrm{Aut}(G)$, and $\mathscr{B}_{U}$ be a basis of $U=\bigcap_{x \in S} L_{x,\Phi}$.\\\indent A map $p:V(G)\to \mathbb{R}^d$ is said to be \emph{$(S,\Phi,\mathscr{B}_{U})$-generic} if the following holds:
If the determinant of any submatrix of $\mathbf{R}(K_{n},p)$ is equal to zero, then the determinant of the
corresponding submatrix of $\mathbf{R}_{\mathscr{B}_{U}}(n,d)$ is (identically) zero.\\\indent The map $p$ is said to be \emph{$(S,\Phi)$-generic} if $p$ is $(S,\Phi,\mathscr{B}_{U})$-generic for some basis $\mathscr{B}_{U}$ of $U$.\\\indent
\indent A framework $(G,p) \in \mathscr{R}_{(G,S,\Phi)}$ is \emph{$(S,\Phi,\mathscr{B}_{U})$-generic} if $p$ is an $(S,\Phi,\mathscr{B}_{U})$-generic map, and $(G,p)$ is \emph{$(S,\Phi)$-generic} if $(G,p)$ is
$(S,\Phi,\mathscr{B}_{U})$-generic for some basis $\mathscr{B}_{U}$ of $U$.}
\end{defin}

\begin{theorem}
\label{indepofbasis} Let $G$ be a graph, $S$ be a symmetry group, and $\Phi$ be a map from $S$ to $\textrm{Aut}(G)$. If $(G,p)\in \mathscr{R}_{(G,S,\Phi)}$ is $(S,\Phi)$-generic, then $(G,p)$ is $(S,\Phi,\mathscr{B}_{U})$-generic for every basis $\mathscr{B}_{U}$ of $U=\bigcap_{x \in S} L_{x,\Phi}$.
\end{theorem}
\textbf{Proof.} Suppose $S$ is a symmetry group in dimension $d$ and the vertex set of $G$ is $V(G)=\{v_{1},v_{2},\ldots,v_{n}\}$. Let $(G,p)\in \mathscr{R}_{(G,S,\Phi)}$ be $(S,\Phi)$-generic, say $(G,p)$ is $(S,\Phi,\mathscr{B}_{U})$-generic, where
$\mathscr{B}_{U}=\{u_{1},\ldots,u_{k}\}$ is a basis of $U$. Let $\mathscr{B^*}_{U}=\{u^*_{1},\ldots,u^*_{k}\}$ be
another basis of $U$. Then we need to show that $(G,p)$ is $(S,\Phi,\mathscr{B^*}_{U})$-generic. Let \begin{displaymath}\left( \begin{array} {c} p_{1} \\ p_{2} \\ \vdots \\ p_{n} \end{array} \right)=t_{1}u_{1}+\ldots+ t_{k}u_{k}=t^*_{1}u^*_{1}+\ldots+ t^*_{k}u^*_{k}\textrm{ ,}\end{displaymath} where $t_{i},t^*_{i}\in \mathbb{R}$ for all $i=1,\ldots,k$. Then there exists an invertible matrix of real numbers $(s_{ij})$ such that \begin{equation}\label{eq:bt}\begin{array} {ccc} t_{1} &=
& s_{11}t^*_{1}+\ldots + s_{1k}t^*_{k}\\ \vdots & \vdots &
\vdots
\\t_{k} &= & s_{k1}t^*_{1}+\ldots + s_{kk}t^*_{k}\textrm{ .}
\end{array}\end{equation}
Let $\mathbf{R}_{\mathscr{B}_{U}}(n,d)$ be the symmetry-adapted indeterminate rigidity matrix corresponding to
$\mathscr{B}_{U}$ with variables $t'_{1},\ldots,t'_{k}$, and $\mathbf{R}_{\mathscr{B^*}_{U}}(n,d)$ be the symmetry-adapted indeterminate rigidity matrix corresponding to $\mathscr{B^*}_{U}$ with variables $t^{*'}_{1},\ldots ,t^{*'}_{k}$.
Then note that if for $i=1,\ldots,k$, we replace the variable $t'_{i}$ in $\mathbf{R}_{\mathscr{B}_{U}}(n,d)$
analogously to \eqref{eq:bt} by \begin{equation}\label{eq:bti} t'_{i} = s_{i1}t^{*'}_{1}+\ldots +
s_{ik}t^{*'}_{k}\textrm{,}\end{equation} then we obtain the matrix
$\mathbf{R}_{\mathscr{B^*}_{U}}(n,d)$.\\\indent If each $t^{*'}_{i}$ in
$\mathbf{R}_{\mathscr{B^*}_{U}}(n,d)$ is replaced by $t^*_{i}$, then, by Remark \ref{symindetrigmatrem}, we obtain the rigidity matrix $\mathbf{R}(K_{n},p)$. Consider the determinant of a submatrix of $\mathbf{R}(K_{n},p)$ which is equal to zero. The determinant of the corresponding submatrix of $\mathbf{R}_{\mathscr{B}_{U}}(n,d)$ is a polynomial in $t'_{1},\ldots,t'_{k}$, say
\begin{equation}\label{eq:pol} \sum a_{(a_{1},\ldots ,a_{k})} t'^{a_{1}}_{1}\cdot \ldots \cdot t'^{a_{k}}_{k} \textrm{, where } a_{(a_{1},\ldots,a_{k})}\in \mathbb{R} \textrm{.}\end{equation} Since $(G,p)$ is $(S,\Phi,\mathscr{B}_{U})$-generic, the polynomial in \eqref{eq:pol} is the zero polynomial. If in \eqref{eq:pol} we replace the variables $t'_{i}$ as in \eqref{eq:bti}, then we again obtain the zero polynomial. On the other hand, this polynomial is the determinant of the corresponding submatrix of
$\mathbf{R}_{\mathscr{B^*}_{U}}(n,d)$. This says that $(G,p)$ is $(S,\Phi,\mathscr{B^*}_{U})$-generic and the proof
is complete. $\square$

Note that it follows directly from Definition \ref{symgeneric} that the set of $(S,\Phi)$-generic realizations of a graph $G$ is an open dense subset of the set $\mathscr{R}_{(G,S,\Phi)}$.\\\indent Moreover, as we will show next, the infinitesimal rigidity properties are the same for all $(S,\Phi)$-generic realizations of $G$.


\begin{lemma}
\label{lemmaforsymgenprop}Let $G$ be a graph with $V(G)=\{v_{1},v_{2},\ldots,v_{n}\}$, $S$ be a symmetry group in dimension $d$, and $\Phi$ be a map from $S$ to $\textrm{Aut}(G)$. If for some framework $(G,p)\in \mathscr{R}_{(G,S,\Phi)}$, the points $p_{1},\ldots,p_{n}$ span an affine subspace of $\mathbb{R}^d$ of dimension $k$, then for any $(S,\Phi)$-generic realization $(G,q)$ of $G$, the points $q_{1},\ldots,q_{n}$ span an affine subspace of $\mathbb{R}^d$ of dimension at least $k$.
\end{lemma}
\textbf{Proof.} Let $(G,p)\in \mathscr{R}_{(G,S,\Phi)}$ be a framework for which the points $p_{1},\ldots,p_{n}$ span an affine subspace of $\mathbb{R}^d$ of dimension $k$. Then there are $k+1$ affinely independent points among $p_{1},\ldots,p_{n}$, say wlog $p_{1},\ldots,p_{k+1}$. Let $A$ be the $k\times d$ matrix defined by
\begin{displaymath}A=\left( \begin{array} {cccc}
(p_{1})_{1}-(p_{2})_{1} & (p_{1})_{2}-(p_{2})_{2} & \ldots & (p_{1})_{d}-(p_{2})_{d}\\
(p_{1})_{1}-(p_{3})_{1} & (p_{1})_{2}-(p_{3})_{2} & \ldots & (p_{1})_{d}-(p_{3})_{d}\\
\vdots & \vdots & \ldots & \vdots\\
(p_{1})_{1}-(p_{k+1})_{1} & (p_{1})_{2}-(p_{k+1})_{2} & \ldots & (p_{1})_{d}-(p_{k+1})_{d}\\
\end{array}
\right)\textrm{.}
\end{displaymath}
Then the rows of $A$ are linearly independent and hence there exists a $k\times k$ submatrix $B$ of $A$ whose determinant is non-zero. Fix a basis $\mathscr{B}_{U}$ of $U=\bigcap_{x \in S} L_{x,\Phi}$ and let $\mathbf{R}_{\mathscr{B}_{U}}(n,d)$ be the symmetry-adapted indeterminate rigidity matrix for $\mathscr{R}_{(G,S,\Phi)}$ corresponding to $\mathscr{B}_{U}$. Then the determinant of the submatrix $B'$ of $\mathbf{R}_{\mathscr{B}_{U}}(n,d)$ which corresponds to $B$ is not identically zero.\\\indent Now, let $(G,q)$ be an $(S,\Phi)$-generic realization of $G$ and suppose the points $q_{1},\ldots,q_{n}$ span an affine subspace of $\mathbb{R}^d$ of dimension $m<k$. Then the matrix $\widehat{A}$ which is obtained from $A$ by replacing each $(p_{i})_{j}$ by $(q_{i})_{j}$ has a non-trivial row dependency, which says that the determinant of every $k\times k$ submatrix of $\widehat{A}$ is equal to zero. This contradicts the fact that $(G,q)$ is $(S,\Phi)$-generic and that the determinant of $B'$ is not identically zero. $\square$

\begin{theorem}
\label{symgenrigthm}
Let $G$ be a graph, $S$ be a symmetry group, and
$\Phi$ be a map from $S$ to $\textrm{Aut}(G)$ such that $\mathscr{R}_{(G,S,\Phi)}\neq \emptyset$. The following are equivalent.
\begin{itemize}
\item[(i)] There exists a framework $(G,p)\in \mathscr{R}_{(G,S,\Phi)}$ that is infinitesimally rigid (independent, isostatic);
\item[(ii)] every $(S,\Phi)$-generic realization of $G$ is infinitesimally rigid (independent, isostatic).
\end{itemize}
\end{theorem}
\textbf{Proof.} Suppose $S$ is a symmetry group in dimension $d$. Let $(G,p) \in \mathscr{R}_{(G,S,\Phi)}$ be
infinitesimally rigid and let $(G,q)$ be an $(S,\Phi)$-generic realization of $G$.\\\indent Suppose first that $|V(G)|\geq d$. Then, by Remark \ref{affine}, the points $p(v)$, $v\in V(G)$, span an affine subspace of $\mathbb{R}^{d}$ of dimension at least $d-1$. Therefore, the trivial infinitesimal motions arising from $d$ translations and $\binom{d}{2}$ rotations of $\mathbb{R}^{d}$ form a basis of the space of trivial infinitesimal motions of $(G,p)$ (see \cite{asiroth, gss} for details), and hence we have
\begin{displaymath} \textrm{rank } \big(\mathbf{R}(G,p)\big)=d |V(G)| -
\binom{d+1}{2}\textrm{.}\end{displaymath}
By the definition of $(S,\Phi)$-generic,
\begin{displaymath} \textrm{rank } \big(\mathbf{R}(G,q)\big) \ge \textrm{rank } \big(\mathbf{R}(G,p)\big).\end{displaymath}

By Lemma \ref{lemmaforsymgenprop}, the points $q(v)$, $v\in V(G)$, also span an affine subspace of $\mathbb{R}^{d}$ of dimension at least $d-1$, which says that $\textrm{nullity }\big(\mathbf{R}(G,q)\big)\geq\binom{d+1}{2}$. Therefore,


\begin{displaymath} \textrm{rank } \big(\mathbf{R}(G,q)\big)\leq d |V(G)| - \binom{d+1}{2}
\textrm{.}\end{displaymath}
It follows that
\begin{displaymath} \textrm{rank } \big(\mathbf{R}(G,q)\big)=d |V(G)| -
\binom{d+1}{2}\textrm{,}\end{displaymath}
and hence $(G,q)$ is infinitesimally rigid.



Suppose now that $|V(G)|\leq d-1$. Then the dimension of the space of trivial infinitesimal motions of $(G,p)$ is strictly smaller than $\binom{d+1}{2}$ (see again \cite{asiroth, gss} for details). Therefore, we have $\textrm{nullity } \big(\mathbf{R}(G,p)\big)<\binom{d+1}{2}$, and hence $\textrm{rank } \big(\mathbf{R}(G,p)\big)>d |V(G)| -\binom{d+1}{2}$. It follows from Theorem \ref{infinrigaff} that $G$ is a complete graph and the points $p(v)$, $v\in V(G)$, are affinely independent. By Lemma \ref{lemmaforsymgenprop}, the points $q(v)$, $v\in V(G)$, must also be affinely independent, and hence $(G,q)$ is infinitesimally rigid.\\\indent
If $(G,p)$ is independent, then it follows from the definition of $(S,\Phi)$-generic that $(G,q)$ is also independent. Therefore, if $(G,p)$ is isostatic, so is $(G,q)$. $\square$

So, being infinitesimally rigid (independent, isostatic) is an $(S,\Phi)$-generic property. This gives rise to

\begin{defin}
\emph{Let $G$ be a graph, $S$ be a symmetry group, and $\Phi$ be a map from $S$ to $\textrm{Aut}(G)$. $G$ is said to be \emph{$(S,\Phi)$-generically infinitesimally rigid (independent, isostatic)} if realizations which are $(S,\Phi)$-generic are infinitesimally rigid (independent, isostatic).}
\end{defin}


Examples \ref{K33ex} and \ref{tpex} show that a graph $G$ which is $(S,\Phi)$-generically isostatic is not necessarily $(S,\Psi)$-generically isostatic, where $\Phi$ and $\Psi$ are two distinct maps from $S$ to $\textrm{Aut}(G)$.\\\indent In Example \ref{K33ex}, $(\mathcal{C}_{s},\Phi_{a})$-generic realizations in $\mathscr{R}_{(K_{3,3},\mathcal{C}_{s},\Phi_{a})}$ are isostatic, whereas all realizations in $\mathscr{R}_{(K_{3,3},\mathcal{C}_{s},\Phi_{b})}$ are not isostatic, because the joints of any realization in $\mathscr{R}_{(K_{3,3},\mathcal{C}_{s},\Phi_{b})}$ lie on a conic section, as we already observed in the beginning of this section.\\\indent In Example \ref{tpex}, the graph $G_{tp}$ is $(\mathcal{C}_{2},\Psi_{a})$-generically isostatic, but none of the realizations in $\mathscr{R}_{(G_{tp},\mathcal{C}_{2},\Psi_{b})}$ are isostatic. This follows from the pure condition for $G_{tp}$, which says that a 2-dimensional realization of $G_{tp}$ is not isostatic if and only if the triangles $p_{1}\, p_{2}\, p_{3}$ and $p_{4}\, p_{5}\, p_{6}$ are perspective from a line \cite{WW}. Equivalently, by Desargues Theorem, a 2-dimensional realization of $G_{tp}$ is not isostatic if and only if the triangles $p_{1}\, p_{2}\, p_{3}$ and $p_{4}\, p_{5}\, p_{6}$ are perspective from a point or at least one of those triangles is degenerate.\\\indent
For an example in 3-space, consider the complete graph $K_{4}$ with $V(K_{4})=\{v_{1},v_{2},v_{3},v_{4}\}$, a symmetry group $\mathcal{C}_{s}=\{Id,s\}$ in dimension 3, and the maps $\Upsilon_{a}$ and $\Upsilon_{b}$ from $\mathcal{C}_{s}$ to $\textrm{Aut}(K_{4})$, where
$\Upsilon_{a}$ maps $Id$ to the identity automorphism $id$ of $K_{4}$ and $s$ to $(v_{1}\, v_{2})(v_{3})(v_{4})$, and $\Upsilon_{b}$ maps both $Id$ and $s$ to $id$. Then $K_{4}$ is $(\mathcal{C}_{s},\Upsilon_{a})$-generically isostatic, but all realizations in $\mathscr{R}_{(K_{4},\mathcal{C}_{s},\Upsilon_{b})}$ are infinitesimally flexible, because all the joints of a realization in $\mathscr{R}_{(K_{4},\mathcal{C}_{s},\Upsilon_{b})}$ must lie in the mirror plane corresponding to $s$ and are therefore coplanar.

\begin{figure}[htp]
\begin{center}
\includegraphics[clip]{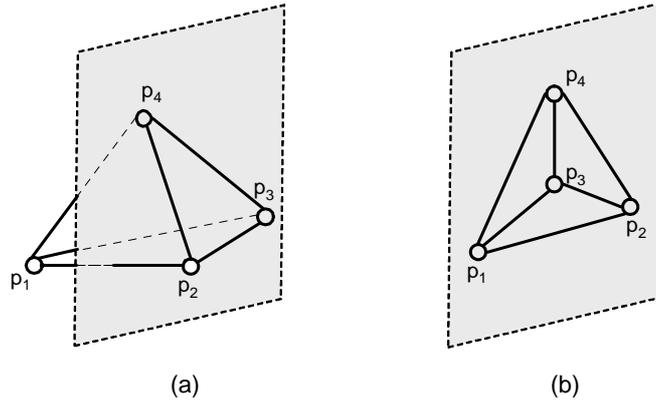}
\end{center}
\caption{\emph{A $3$-dimensional realization of $(K_{4},\mathcal{C}_{s})$ of type $\Upsilon_{a}$ (a) and of type $\Upsilon_{b}$ (b).}}
\label{tettypes}
\end{figure}

\begin{remark}
\label{orbits}
\emph{Let $G$ be a graph, $S$ be a symmetry group in dimension $d$, and $\Phi:S\to \textrm{Aut}(G)$ be a homomorphism. Frameworks in the set  $\mathscr{R}_{(G,S,\Phi)}$, particularly $(S,\Phi)$-generic realizations of $G$, can then be visualized in a very intuitive way via the following approach.\\\indent The map of $S\times V(G)$ onto $V(G)$ that sends $(x,v)$ to $\Phi(x)(v)$ defines a group action on $V(G)$ and the orbits $S v=\{\Phi(x)(v)|\, x\in S\}$ form a partition of $V(G)$. Let $\{v_{1},\ldots,v_{r}\}$ be a subset of $V(G)$ obtained by choosing one representative from each of these orbits, and recall from Definition \ref{symel} that for every $x\in S$, $F_{x}$ denotes the symmetry element corresponding to $x$. If for $i=1,\ldots,r$, we define \begin{displaymath} F(v_{i})=\bigcap_{x\in S\textrm{ with }\Phi(x)(v_{i})=v_{i}}F_{x}\textrm{,}\end{displaymath} then for every framework $(G,p)\in \mathscr{R}_{(G,S,\Phi)}$ and for every $i\in\{1,\ldots,r\}$, the point $p(v_{i})$ must be contained in the subspace $F(v_{i})$ of $\mathbb{R}^{d}$. \\\indent Note that the positions of \emph{all} joints of a framework $(G,p)\in \mathscr{R}_{(G,S,\Phi)}$ are uniquely determined by the positions $p(v_{1}),\ldots,p(v_{r})$  of the joints $\big(v_{1},p(v_{1})\big),\ldots,\big(v_{r},p(v_{r})\big)$ and the symmetry constraints imposed by $S$ and $\Phi$. In other words, we may construct frameworks in $\mathscr{R}_{(G,S,\Phi)}$ by first choosing a point $p(v_{i})\in F(v_{i})$ for each $i=1,\ldots, r$ and then letting $S$ and $\Phi$ determine the positions of the remaining joints. In particular, note that we obtain an $(S,\Phi)$-generic framework $(G,p)$ in this way for almost all choices of points $p(v_{i})$ that satisfy $p(v_{i})\in F(v_{i})$ for $i=1,\ldots, r$.\\\indent Consider, for example, the set $\mathscr{R}_{(K_{4},\mathcal{C}_{s},\Upsilon_{a})}$ an element of which is shown in Figure \ref{tettypes} (a). The orbits for the group action from $\mathcal{C}_{s}\times V(K_{4})$ onto $V(K_{4})$ are given by $\{v_{1},v_{2}\}$, $\{v_{3}\}$, and $\{v_{4}\}$. If $(K_{4},p)$ is a framework in $\mathscr{R}_{(K_{4},\mathcal{C}_{s},\Upsilon_{a})}$, then both $p_{3}$ and $p_{4}$ must be contained in the mirror plane $F_{s}$ of $s$, because $F(v_{3})=F(v_{4})=F_{Id}\cap F_{s}=F_{s}$. Furthermore, since $v_{1}$ and $v_{2}$ belong to the same orbit, the position of the point $p_{2}$ is uniquely determined by the position of $p_{1}$ and the symmetry constraints imposed by $\mathcal{C}_{s}$ and $\Upsilon_{a}$. Since $F(v_{1})=F_{Id}=\mathbb{R}^3$, the point $p_{1}$ may be chosen to be any point in $\mathbb{R}^3$; however, if $p_{1}$ lies in the mirror plane of $s$, then $p_{1}=p_{2}$, in which case $(K_{4},p)$ is not a framework.\\\indent For a further elaboration on this approach we refer the reader to \cite{BS2}. }
\end{remark}


We conclude this section by giving a few more interesting properties of $(S,\Phi)$-generic frameworks.

\begin{theorem}
\label{symgenprop} Let $G$ be a graph, $S$ be a symmetry group, and $\Phi$ be a map from $S$ to $\textrm{Aut}(G)$. Further, let $(G,p)\in \mathscr{R}_{(G,S,\Phi)}$, $S'$ be a subgroup of $S$, and $\Phi'=\Phi|_{S'}$. If $(G,p)$ is $(S',\Phi')$-generic, then $(G,p)$ is also $(S,\Phi)$-generic.
\end{theorem}
\textbf{Proof.} Suppose $S$ is a symmetry group in dimension $d$ and $G$ is a graph with $n$ vertices. Let $(G,p)\in \mathscr{R}_{(G,S,\Phi)}$ be $(S',\Phi')$-generic. We fix a basis $\mathscr{B}_{U}=\{u_{1},\ldots,u_{k}\}$ of $U = \bigcap_{x \in S}L_{x,\Phi}$ and extend it to a basis $\mathscr{B}_{U'}=\{u_{1},\ldots,u_{k},u_{k+1},\ldots,u_{m}\}$ of
$U' = \bigcap_{x \in S'} L_{x,\Phi'}$. Consider a submatrix of the rigidity matrix $\mathbf{R}(K_{n},p)$ whose determinant is equal to zero. We need to show that the determinant $\Delta$ of the corresponding submatrix of the symmetry-adapted indeterminate rigidity matrix
$\mathbf{R}_{\mathscr{B}_{U}}(n,d)$ is identically zero.
\\ \indent Let $\Delta'=\sum a_{(a_{1},\ldots ,a_{m})}t'^{a_{1}}_{1}\cdot \ldots \cdot t'^{a_{m}}_{m} \textrm{ , where }
a_{(a_{1},\ldots,a_{m})}\in \mathbb{R}$, be the determinant of the corresponding submatrix of
$\mathbf{R}_{\mathscr{B}_{U'}}(n,d)$. Then $\Delta'$ is the zero polynomial since $(G,p)$ is $(S',\Phi')$-generic. But note that $\Delta$ is a polynomial that is obtained from $\Delta'$ by deleting all those terms in $\Delta'$ that have one or
more variables in $\{t'_{k+1},\ldots,t'_{m}\}$. Thus, $\Delta$ is also the zero polynomial. $\square$

The converse of Theorem \ref{symgenprop} does not hold.

\begin{figure}[htp]
\begin{center}
\includegraphics[clip]{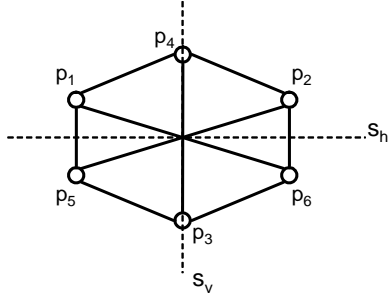}
\end{center}
\caption{\emph{A realization of $K_{3,3}$ that is $(\mathcal{C}_{2v},\Phi)$-generic, but not $(\mathcal{C}_{s},\Phi_{a})$-generic, where $\mathcal{C}_{s}$ is the subgroup of $\mathcal{C}_{2v}$ generated by $s_{v}$ and $\Phi_{a}=\Phi|_{\mathcal{C}_{s}}$.}}
\label{gs}
\end{figure}

\begin{examp}
\label{highersymgen}
\emph{The realization $(K_{3,3},p)$ in Figure \ref{gs} is $(\mathcal{C}_{2v},\Phi)$-generic, where $\mathcal{C}_{2v}=\{Id, C_{2}, s_{h},s_{v}\}$ is a symmetry group in dimension 2 and $\Phi: \mathcal{C}_{2v}\to \textrm{Aut}(K_{3,3})$ is defined by
\begin{eqnarray} \Phi(Id) & = & id\nonumber\\ \Phi(C_{2}) & = & (v_{1}\,v_{6})(v_{2}\, v_{5})(v_{3}\,v_{4})\nonumber\\ \Phi(s_{h}) & = & (v_{1}\,v_{5})(v_{2}\, v_{6})(v_{3}\,v_{4})\nonumber\\ \Phi(s_{v}) & = & (v_{1}\,v_{2})(v_{5}\,v_{6})(v_{3}) (v_{4})
\nonumber\textrm{.}\end{eqnarray}
However, $(K_{3,3},p)$ is not $(\mathcal{C}_{s},\Phi_{a})$-generic, where $\mathcal{C}_{s}$ is the subgroup of $\mathcal{C}_{2v}$ generated by $s_{v}$ and $\Phi_{a}=\Phi|_{\mathcal{C}_{s}}$ is the map we defined in Example \ref{K33ex}.}
\end{examp}

\begin{cor}
\label{genericproperties}
Let $G$ be a graph, $S$ be a symmetry group, and $\Phi$ be a map from $S$ to $\textrm{Aut}(G)$. If $(G,p)\in\mathscr{R}_{(G,S,\Phi)}$ is generic (in the sense of Definition \ref{generic}), then $(G,p)$ is also $(S,\Phi)$-generic.
\end{cor}
\textbf{Proof.} Suppose $S$ is a symmetry group in dimension $d$ and $G$ is a graph with $n$ vertices. Let $(G,p)\in\mathscr{R}_{(G,S,\Phi)}$ be generic. Then $\Phi$ maps the symmetry operation $Id\in S$ to the identity automorphism $id$ of $G$, for otherwise the map $q$ of every realization $(G,q)$ in $\mathscr{R}_{(G,S,\Phi)}$ is non-injective, contradicting the fact that $(G,p)\in\mathscr{R}_{(G,S,\Phi)}$ is generic. So, $\Phi|_{\mathcal{C}_{1}}=I$, where $I:\mathcal{C}_{1}\to \textrm{Aut}(G)$ maps $Id$ to $id$, and we have $\bigcap_{x \in \mathcal{C}_{1}}L_{x,I}=L_{Id,I}=\mathbb{R}^{dn}$.\\\indent Now, observe that the indeterminate rigidity matrix $\mathbf{R}(n,d)$ is equal to the symmetry-adapted indeterminate rigidity matrix $\mathbf{R}_{\mathscr{B}_{\mathbb{R}^{dn}}}(n,d)$, where $\mathscr{B}_{\mathbb{R}^{dn}}$ is the canonical basis of $\mathbb{R}^{dn}$. Therefore, $(G,p)$ is generic if and only if $(G,p)$ is $(\mathcal{C}_{1},I)$-generic.\\\indent The result now follows immediately from Theorem \ref{symgenprop}. $\square$

The converse of Corollary \ref{genericproperties} is of course false. A $(\mathcal{C}_{s},\Phi_{b})$-generic realization of $K_{3,3}$, for example, where $\mathcal{C}_{s}$ and $\Phi_{b}$ are as in Example \ref{K33ex}, has all of its joints on a conic section and is therefore not generic.

\section{Of what types $\Phi$ can a framework be?}

Let $(G,p)$ be a $d$-dimensional framework with point group symmetry $P$. Then $(G,p)\in \mathscr{R}_{(G,S)}$ for every subgroup $S$ of $P$.\\\indent Fix a subgroup $S$ of $P$. Then it follows from Theorem \ref{symcharac} that there exists a map $\Phi:S \to \textrm{Aut}(G)$ such that $(G,p)\in \mathscr{R}_{(G,S,\Phi)}$. The following examples show that it is possible for $(G,p)\in \mathscr{R}_{(G,S)}$ to be of more than just one such type $\Phi$. Note that each of those examples is a non-injective realization.


\begin{figure}[htp]
\begin{center}
\includegraphics[clip]{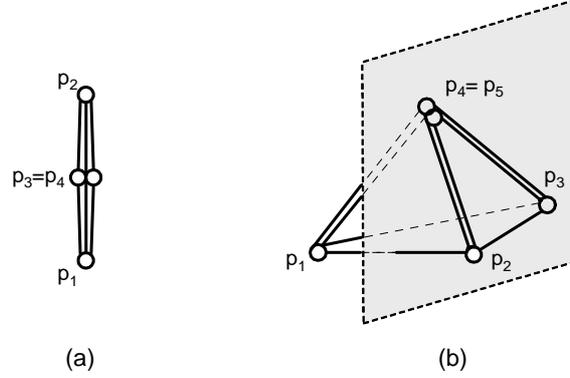}
\end{center}
\caption{\emph{A realization of $(G_{t},\mathcal{C}_{2})$ of type $\Theta_{a}$ and $\Theta_{b}$ (a) and a realization of $(G_{bp},\mathcal{C}_{s})$ of type $\Xi_{a}$ and $\Xi_{b}$ (b).}}
\label{noninjreal}
\end{figure}

\begin{examp}
\label{tsharingedge}\emph{Let $G_{t}$ be the graph of two triangles sharing an edge and $\mathcal{C}_{2}=\{Id,C_{2}\}$ be the half-turn symmetry group in dimension 2. Figure \ref{noninjreal} (a) shows a realization $(G_{t},p)$ of $(G_{t},\mathcal{C}_{2})$ that is of type $\Theta_{a}$ as well as $\Theta_{b}$, where $\Theta_{a}:\mathcal{C}_{2}\to \textrm{Aut}(G_{t})$ is defined by \begin{eqnarray} \Theta_{a}(Id) & = & id\nonumber\\ \Theta_{a}(C_{2}) & = & (v_{1}\,v_{2})(v_{3})(v_{4})\textrm{, }\nonumber
\end{eqnarray}
and $\Theta_{b}:\mathcal{C}_{2}\to \textrm{Aut}(G_{t})$ is defined by \begin{eqnarray} \Theta_{b}(Id) & = & id\nonumber\\ \Theta_{b}(C_{2}) & = & (v_{1}\,v_{2})(v_{3}\,v_{4})\textrm{.}\nonumber
\end{eqnarray}}
\end{examp}

\begin{examp}
\label{bipyrnoninj}\emph{ Consider the graph $G_{bp}$ of a triangular bipyramid and a symmetry group $\mathcal{C}_{s}=\{Id,s\}$ in dimension 3. The framework $(G_{bp},p)$ in Figure \ref{noninjreal} (b) is a realization of $(G_{bp},\mathcal{C}_{s})$ that is of type $\Xi_{a}$ as well as $\Xi_{b}$, where $\Xi_{a}:\mathcal{C}_{s}\to \textrm{Aut}(G_{bp})$ and $\Xi_{b}:\mathcal{C}_{s}\to \textrm{Aut}(G_{bp})$ are defined as in Example \ref{3Dtypes}.}
\end{examp}

Since for a given framework $(G,p)$ in a set of the form $\mathscr{R}_{(G,S)}$, the specification of a type $\Phi:S\to \textrm{Aut}(G)$ plays a key role in a symmetry-based rigidity analysis of $(G,p)$, it is natural to ask how we can find all the types $\Phi$ of $(G,p)$, how these types are related to each other and under what conditions $(G,p)$ is of a unique type.


The following definition is essential to answer all of these questions.

\begin{defin}
\emph{Let $(G,p)$ be a framework. Then we denote \emph{$\textrm{Aut}(G,p)$} to be the set of all $\alpha\in \textrm{Aut}(G)$ which satisfy $p(v)=p\big(\alpha(v)\big)$ for all $v\in V(G)$.}
\end{defin}

Given a framework $(G,p)$ and an automorphism $\alpha\in \textrm{Aut}(G,p)$, it is easy to see that only vertices of $G$ that have the same image under $p$ can possibly belong to the same permutation cycle of $\alpha$. In particular, for every framework $(G,p)$ with an injective map $p$, we have $\textrm{Aut}(G,p)=\{id\}$, as we will see in the proof of Corollary \ref{autgpinjective}.\\\indent For the framework $(G_{t},p)$ in Figure \ref{noninjreal} (a), we have $\textrm{Aut}(G_{t},p)=\{id, (v_{3}\, v_{4})(v_{1})(v_{2})\}$ and for the framework $(G_{bp},p)$ in Figure \ref{noninjreal} (b), we have $\textrm{Aut}(G_{bp},p)=\{id, (v_{4}\, v_{5})(v_{1})(v_{2})(v_{3})\}$. \\\indent
Clearly, $\textrm{Aut}(G,p)$ is a subgroup of $\textrm{Aut}(G)$.

\begin{theorem}
\label{autgp1}
Let $G$ be a graph, $S$ be a symmetry group, and $\Phi$ be a map from $S$ to $\textrm{Aut}(G)$. Further, let $(G,p)\in \mathscr{R}_{(G,S,\Phi)}$ and $x\in S$. Then $\Phi(x)\textrm{Aut}(G,p)=\textrm{Aut}(G,p)\Phi(x)$, and an automorphism $\alpha$ of $G$ satisfies $x\big(p(v)\big)=p\big(\alpha(v)\big)$ for all $v\in V(G)$ if and only if $\alpha$ is an element of $\Phi(x)\textrm{Aut}(G,p)$.
\end{theorem}
\textbf{Proof.} First, we show that $\Phi(x)\textrm{Aut}(G,p)=\textrm{Aut}(G,p)\Phi(x)$. Since the cosets $\Phi(x)\textrm{Aut}(G,p)$ and $\textrm{Aut}(G,p)\Phi(x)$ have the same cardinality, it suffices to show that $\Phi(x)\textrm{Aut}(G,p)\subseteq \textrm{Aut}(G,p)\Phi(x)$. Let $\alpha\in \Phi(x)\textrm{Aut}(G,p)$, say $\alpha=\Phi(x)\circ \beta$, where $\beta\in \textrm{Aut}(G,p)$. Then for $v\in V(G)$, we have \begin{displaymath}x\big(p(v)\big)=x\Big(p\big(\beta(v)\big)\Big)=p\Big(\Phi(x)\big(\beta(v)\big)\Big)=p\big(\alpha(v)\big)
\textrm{.}\end{displaymath} Since we also have $x\big(p(v)\big)=p\big(\Phi(x)(v)\big)$, it follows that $p\big(\alpha(v)\big)=p\big(\Phi(x)(v)\big)$ for all $v\in V(G)$. Therefore, \begin{displaymath}p\Big(\alpha\circ \big(\Phi(x)\big)^{-1}(v)\Big)=p(v) \textrm{ for all } v\in V(G)\textrm{,}\end{displaymath} and hence $\alpha \circ (\Phi(x))^{-1}\in \textrm{Aut}(G,p)$. Thus, $\alpha\in \textrm{Aut}(G,p)\Phi(x)$.\\\indent Now, $\alpha\in \textrm{Aut}(G)$ satisfies \begin{displaymath}x\big(p(v)\big)=p\big(\alpha(v)\big) \textrm{ for all } v\in V(G)\end{displaymath} if and only if \begin{displaymath}p\big(\alpha(v)\big)=p\big(\Phi(x)(v)\big)\textrm{ for all } v\in V(G)\end{displaymath} if and only if \begin{displaymath}p\Big(\alpha\circ \big(\Phi(x)\big)^{-1}(v)\Big)=p(v)\textrm{ for all } v\in V(G)\end{displaymath} if and only if \begin{displaymath}\alpha \circ \big(\Phi(x)\big)^{-1}\in \textrm{Aut}(G,p)\end{displaymath} if and only if \begin{displaymath}\alpha\in \textrm{Aut}(G,p)\Phi(x)=\Phi(x)\textrm{Aut}(G,p)\textrm{.} \qquad\square\end{displaymath}

\begin{cor}
\label{autgpid}
Let $G$ be a graph, $S$ be a symmetry group, $\Phi$ be a map from $S$ to $\textrm{Aut}(G)$, and $(G,p)\in \mathscr{R}_{(G,S,\Phi)}$. Then for every $\Psi:S \to \textrm{Aut}(G)$ distinct from $\Phi$, we have $(G,p)\notin \mathscr{R}_{(G,S,\Psi)}$ if and only if $\textrm{Aut}(G,p)=\{id\}$.
\end{cor}
\textbf{Proof.} It follows directly from Theorem \ref{autgp1} that $\textrm{Aut}(G,p)=\{id\}$ if and only if for every $x\in S$, the automorphism $\Phi(x)$ is the only automorphism of $G$ that satisfies $x\big(p(v)\big)=p\big(\Phi(x)(v)\big)$ for all $v\in V(G)$. $\square$

Corollary \ref{autgpid} asserts that the type $\Phi:S\to \textrm{Aut}(G)$ of a framework $(G,p)\in \mathscr{R}_{(G,S)}$ is unique if and only if $\textrm{Aut}(G,p)$ only contains the identity automorphism of $G$. In particular, we have the following result.

\begin{cor}
\label{autgpinjective}
Let $G$ be a graph, $S$ be a symmetry group, and $\Phi$ be a map from $S$ to $\textrm{Aut}(G)$. If the map $p$ of a framework $(G,p)\in\mathscr{R}_{(G,S,\Phi)}$ is injective, then $(G,p)\notin \mathscr{R}_{(G,S,\Psi)}$ for every $\Psi:S \to \textrm{Aut}(G)$ distinct from $\Phi$.
\end{cor}
\textbf{Proof.} Let $\alpha$ be an element of $\textrm{Aut}(G,p)$. Then we have $p(v)=p\big(\alpha(v)\big)$ for all $v\in V(G)$, and since $p$ is injective it follows that $v=\alpha(v)$ for all $v\in V(G)$. Thus, $\alpha$ is the identity automorphism of $G$ and the result follows from Corollary \ref{autgpid}.
$\square$

The following examples show that the converse of Corollary \ref{autgpinjective} does not hold, that is, a framework $(G,p)\in \mathscr{R}_{(G,S)}$ that is of a unique type $\Phi:S\to \textrm{Aut}(G)$ can possibly have a non-injective map $p$.

\begin{figure}[htp]
\begin{center}
\includegraphics[clip]{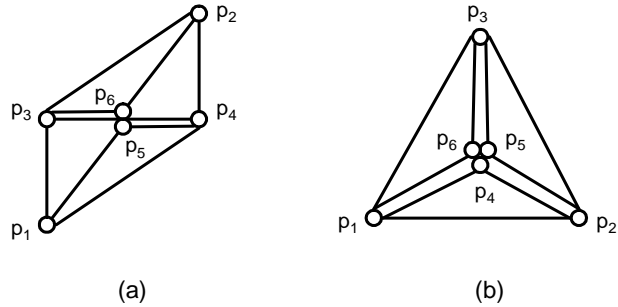}
\end{center}
\caption{\emph{Non-injective realizations $(G,p)$ with $\textrm{Aut}(G,p)=\{id\}$.}}
\label{ninid}
\end{figure}

\begin{examp}
\label{noninjuniqetype}\emph{The framework $(G,p)$ in Figure \ref{ninid} (a) is a non-injective realization of $(G,\mathcal{C}_{2})$ (since $p_{5}=p_{6}$) with $\textrm{Aut}(G,p)=\{id\}$. So, $(G,p)\in \mathscr{R}_{(G,\mathcal{C}_{2})}$ is of the unique type $\Phi:\mathcal{C}_{2}\to \textrm{Aut}(G)$, where $\Phi(Id)=id$ and $\Phi(C_{2})=(v_{1}\, v_{2})(v_{3}\, v_{4})(v_{5}\, v_{6})$.}
\end{examp}

\begin{examp}
\label{noninjuniqetype2}\emph{The framework $(G,p)$ in Figure \ref{ninid} (b) is a non-injective realization of $(G,\mathcal{C}_{3})$ (since $p_{4}=p_{5}=p_{6}$) with $\textrm{Aut}(G,p)=\{id\}$. So, $(G,p)\in \mathscr{R}_{(G,\mathcal{C}_{3})}$ is of the unique type $\Phi:\mathcal{C}_{3}\to \textrm{Aut}(G)$, where $\Phi$ is the homomorphism defined by $\Phi(C_{3})=(v_{1}\, v_{2}\,v_{3})(v_{4} \, v_{5}\, v_{6})$.}
\end{examp}

\begin{remark}
\label{uniquetypeisgeneric}
\emph{Let $(G,p)\in \mathscr{R}_{(G,S,\Phi)}$ be a framework with $\textrm{Aut}(G,p)=\{id\}$ and let $(G,q)\in \mathscr{R}_{(G,S,\Phi)}$ be an $(S,\Phi)$-generic framework. It follows immediately from the definition of $(S,\Phi)$-generic (Definition \ref{symgeneric}) that two joints $(v_{i},q_{i})$ and $(v_{j},q_{j})$ of $(G,q)$ can only satisfy $q_{i}=q_{j}$ if $p_{i}=p_{j}$. This says that $(G,q)$ also satisfies $\textrm{Aut}(G,q)=\{id\}$. Therefore, by Corollary \ref{autgpid}, being of a unique type is an $(S,\Phi)$-generic property.}
\end{remark}


\begin{remark}
\label{diftypesgen}
\emph{If a framework $(G,p)\in \mathscr{R}_{(G,S)}$ is of distinct types $\Phi_{1}, \ldots \Phi_{k}$, where $k\geq 2$, then  $(G,p)$ is not $(S,\Phi_{t})$-generic for some $t\in\{1,\ldots, k\}$, as the following argument shows.\\\indent Suppose to the contrary that $(G,p)$ is $(S,\Phi_{i})$-generic for all $i=1,\ldots,k$ and let $l\in\{1,\ldots,k\}$. Since $\textrm{Aut}(G,p)\neq\{id\}$, there exist vertices $v\neq w$ of $G$ such that $p(v)=p(w)$ and $\alpha(v)=w$ for some $\alpha\in \textrm{Aut}(G,p)$. Since $(G,p)$ is $(S,\Phi_{l})$-generic, there must exist non-trivial symmetry operations $x,y\in S$ such that $\Phi_{l}(x)(v)=v$ and $\Phi_{l}(y)(w)=w$ (and the symmetry elements corresponding to $x$ and $y$ must be the origin $0=p(v)=p(w)$). If for each $x\in S$ with $\Phi_{l}(x)(v)=v$, we replace $\Phi_{l}(x)$ by $\alpha\circ\Phi_{l}(x)$, then we obtain a map $\Phi_{t}$, $t\neq l$, with the property that for all $x\in S$, $\Phi_{t}(x)(v)\neq v$. Thus, $(G,p)$ is not $(S,\Phi_{t})$-generic, a contradiction.\\\indent As an example, consider the framework $(G_{t},p)$ in Figure \ref{noninjreal} (a). $(G_{t},p)$ is $(\mathcal{C}_{2},\Theta_{a})$-generic, but not $(\mathcal{C}_{2},\Theta_{b})$-generic, because $p_{3}=p_{4}$ and $\Theta_{b}(v_{3})=v_{4}$ (see Example \ref{tsharingedge}).\\\indent The framework in Figure \ref{noninjreal} (b) is a realization of $(G_{bp},\mathcal{C}_{s})$ of type $\Xi_{a}$ and $\Xi_{b}$ which is neither $(\mathcal{C}_{s},\Xi_{a})$-generic nor $(\mathcal{C}_{s},\Xi_{b})$-generic, because $p_{4}=p_{5}$ (see Example \ref{bipyrnoninj}).}
\end{remark}

\section{When is a type $\Phi$ of a framework a homomorphism?}

As we mentioned in the introduction, many important results concerning the rigidity properties of symmetric frameworks can be obtained by using group representation theory techniques. Given a set $\mathscr{R}_{(G,S,\Phi)}$ of symmetric realizations, two particular matrix representations of the group $S$, the so-called internal and external representation (see \cite{FGsymmax, KG2, KG3,BS4, BS1}), are basic to these techniques. Both of these matrix representations depend on the graph $G$ and the map $\Phi:S\to \textrm{Aut}(G)$; however, they are only matrix representations of $S$, if $\Phi$ is a homomorphism \cite{BS1, BS4}. So, only in the case where $\Phi$ is a homomorphism we can use group representation theory techniques to further analyze the rigidity properties of a framework in $\mathscr{R}_{(G,S,\Phi)}$. Therefore, we now turn our attention to the important question under what conditions a type $\Phi$ of a given framework in $\mathscr{R}_{(G,S)}$ is in fact a group homomorphism (rather than just a map).


\begin{theorem}
\label{hom}
Let $S$ be a symmetry group and $(G,p)$ be a framework in $\mathscr{R}_{(G,S)}$ with $\textrm{Aut}(G,p)=\{id\}$. Then the unique map $\Phi:S\to \textrm{Aut}(G)$ for which $(G,p)\in \mathscr{R}_{(G,S,\Phi)}$ is a homomorphism.
\end{theorem}
\textbf{Proof.} Let $x$ and $y$ be any two elements of $S$. Then $\Phi(y)\circ\Phi(x)\in \textrm{Aut}(G)$ satisfies \begin{displaymath}\big(y\circ x\big) \big(p(v)\big)=y\Big(p\big(\Phi(x)(v)\big)\Big)=p\Big(\big(\Phi(y)\circ\Phi(x)\big)(v)\Big) \textrm{ for all } v\in V(G)\end{displaymath} and, by Corollary \ref{autgpid}, $\Phi(y)\circ\Phi(x)$ is the only automorphism of $G$ with this property. Thus, $\Phi(y\circ x)=\Phi(y)\circ\Phi(x)$. $\square$


In particular, it follows from Corollary \ref{autgpinjective} and Theorem \ref{hom} that if the map $p$ of $(G,p) \in \mathscr{R}_{(G,S)}$ is injective, then the unique type $\Phi$ of $(G,p)$ is a group homomorphism.

\begin{theorem}
\label{sufhom}
Let $S$ be a symmetry group, $\Phi:S\to \textrm{Aut}(G)$ be a map, and $(G,p)$ be a framework in $\mathscr{R}_{(G,S,\Phi)}$.
\begin{itemize}
\item[(i)] If $\Phi$ is a homomorphism, then $\Phi(S)$ is a subgroup of $\textrm{Aut}(G)$;
\item[(ii)] if $\Phi(S)$ is a subgroup of $\textrm{Aut}(G)$ and $\Phi(x)=\Phi(y)$ whenever $\Phi(y)\in \Phi(x)\textrm{Aut}(G,p)$, then $\Phi$ is a homomorphism.
\end{itemize}
\end{theorem}
\textbf{Proof.} $(i)$ It is a standard result in algebra that the homomorphic image of a group is again a group.\\\indent $(ii)$ Let $x$ and $y$ be any two elements of $S$. By the same argument as in the proof of Theorem \ref{hom}, we have \begin{displaymath}\big(y\circ x\big)\big(p(v)\big)=p\Big(\big(\Phi(y)\circ\Phi(x)\big)(v)\Big)\textrm{ for all } v\in V(G)\textrm{.}\end{displaymath} It follows from Theorem \ref{autgp1} that $\Phi(y\circ x)\in \big(\Phi(y)\circ\Phi(x)\big)\textrm{Aut}(G,p)$. By assumption, $\Phi(S)$ contains at most one element of each of the cosets of $\textrm{Aut}(G,p)$. Since $\Phi(S)$ is a group, the element of the coset $\big(\Phi(y)\circ\Phi(x)\big)\textrm{Aut}(G,p)$ that lies in $\Phi(S)$ must be $\Phi(y)\circ\Phi(x)$. It follows that $\Phi(y\circ x)=\Phi(y)\circ\Phi(x)$ and the proof is complete. $\square$

For a framework $(G,p) \in \mathscr{R}_{(G,S)}$ with $\textrm{Aut}(G,p)\neq \{id\}$, there does not necessarily exist any
homomorphism $\Phi:S\to \textrm{Aut}(G)$ for which $(G,p) \in\mathscr{R}_{(G,S,\Phi)}$, as the following examples illustrate.

\begin{figure}[htp]
\begin{center}
\includegraphics[clip]{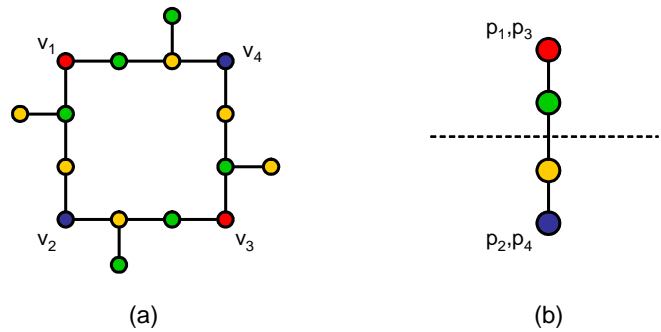}
\end{center}
\caption{\emph{A graph $G$ (a) and a realization  $(G,p)\in\mathscr{R}_{(G,\mathcal{C}_{s})}$ (b) for which there does not exist a homomorphism $\Phi:\mathcal{C}_{s}\to \textrm{Aut}(G)$ so that $(G,p)$ is of type $\Phi$.}}
\label{hom1}
\end{figure}

\begin{examp}
\label{homex1}
\emph{Consider the graph $G$ and the 2-dimensional realization $(G,p)$ of $G$ shown in Figure \ref{hom1} (a) and (b), respectively. Let $s$ be the reflection whose mirror line is shown in Figure \ref{hom1} (b). All vertices of $G$ that are illustrated with the same color in Figure \ref{hom1} have the same image under $p$. Observe that the `$\frac{1}{4}$-turn-automorphism' $\sigma$ of $G$ that permutes the vertices $v_{1},v_{2},v_{3}$ and $v_{4}$ according to the cycle $(v_{1}\,v_{2}\,v_{3}\,v_{4})$ satisfies $s\big(p(v_{i})\big)=p\big(\sigma(v_{i})\big)$ for all $v_{i}\in V(G)$. Thus, $s$ is a symmetry operation of $(G,p)$, and hence
$(G,p)$ is an element of $\mathscr{R}_{(G,\mathcal{C}_{s})}$, where $\mathcal{C}_{s}=\{Id,s\}$.\\ \indent Note that
$\textrm{Aut}(G,p)=\{id,\sigma^2\}$. Therefore, by Theorem \ref{autgp1}, $id$ and $\sigma$ are the two automorphisms of $G$ that can turn $Id\in\mathcal{C}_{s}$ into a symmetry operation of $(G,p)$. Similarly, either one of the elements of $\sigma \textrm{Aut}(G,p)=\{\sigma,\sigma^3\}$ can turn $s\in\mathcal{C}_{s}$ into a symmetry operation of $(G,p)$. It now follows from Theorem \ref{sufhom} $(i)$ that there does not exist any homomorphism $\Phi:\mathcal{C}_{s}\to \textrm{Aut}(G)$ such that $(G,p)\in \mathscr{R}_{(G,\mathcal{C}_{s})}$ is of type $\Phi$, because we cannot choose two elements, one from each of the cosets $\textrm{Aut}(G,p)$ and $\sigma \textrm{Aut}(G,p)$, that form a subgroup of $\textrm{Aut}(G)$.}
\end{examp}


\begin{figure}[htp]
\begin{center}
\includegraphics[clip]{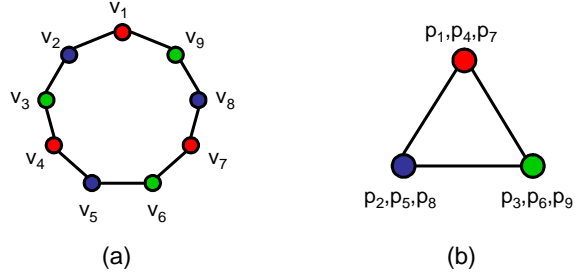}
\end{center}
\caption{\emph{A graph $G$ (a) and a realization  $(G,p)\in\mathscr{R}_{(G,\mathcal{C}_{3})}$ (b) for which there does not exist a homomorphism $\Phi:\mathcal{C}_{3}\to \textrm{Aut}(G)$ so that $(G,p)$ is of type $\Phi$.}}
\label{hom2}
\end{figure}

\begin{examp}
\label{homex2}
\emph{Consider the graph $G$ and the 2-dimensional realization $(G,p)$ of $G$ shown in Figure \ref{hom2} (a) and (b), respectively. As in the previous example, all vertices of $G$ that are illustrated with the same color in Figure \ref{hom2} have the same image under $p$. Note that $(G,p)$ is an element of $\mathscr{R}_{(G,\mathcal{C}_{3})}$, where $\mathcal{C}_{3}=\{Id,C_{3},C_{3}^2\}$ is a symmetry group in dimension 2, because the automorphism $\gamma=(v_{1}\,v_{2}\,\ldots\, v_{9})$ of
$G$ satisfies $C_{3}\big(p(v_{i})\big)=p\big(\gamma(v_{i})\big)$ for all $v_{i}\in V(G)$ and
$\gamma^2$ satisfies $C_{3}^2\big(p(v_{i})\big)=p\big(\gamma^2(v_{i})\big)$ for all $v_{i}\in V(G)$.\\\indent We have $\textrm{Aut}(G,p)=\{id,\gamma^3, \gamma^6\}$, and hence $\gamma \textrm{Aut}(G,p)=\{\gamma,\gamma^4, \gamma^7\}$ and $\gamma^2
\textrm{Aut}(G,p)=\{\gamma^2,\gamma^5, \gamma^8\}$. Since $C_{3}\in\mathcal{C}_{3}$ has order 3 and each element in $\gamma
\textrm{Aut}(G,p)$ has order 9 it follows that there does not exist any homomorphism $\Phi:\mathcal{C}_{3}\to \textrm{Aut}(G)$ such that
$(G,p)$ is of type $\Phi$.}
\end{examp}

It follows that for the frameworks in the above examples, we cannot use group representation theory techniques to analyze their rigidity properties. Note that Examples \ref{homex1} and \ref{homex2} can easily be extended to obtain further examples of frameworks $(G,p)$ and symmetry groups $S$ with the property that there exists no homomorphism $\Phi:S\to \textrm{Aut}(G)$ for which $(G,p) \in\mathscr{R}_{(G,S,\Phi)}$.

\section{Extensions}

While we have only considered \emph{infinitesimal} rigidity properties of symmetric frameworks in this paper, our results can of course also be used for a symmetry-based attack on problems in static rigidity (which is the equivalent dual of infinitesimal rigidity) or rigidity (i.e., the theory of finite flexes of frameworks), as demonstrated in \cite{BS4, BS1, BS2}, for example.\\\indent
Also, while we have restricted our attention to symmetric \emph{bar and joint frameworks}, i.e., realizations of graphs in Euclidean $d$-space, all the definitions and results of this paper can easily be extended to more general classes of symmetric structures.\\\indent For example, we may replace the underlying combinatorial structure of a framework, i.e., a graph, by a more general structure, such as a multigraph or a hypergraph. For any such generalization of a framework, we may then define a symmetry operation, symmetry element, and symmetry group in the analogous way as in Section 2, and thus also establish an analogous classification of these structures. In particular, we obtain linear equations similar to the ones in Definition \ref{classif} in this way so that we can also introduce a symmetry-adapted notion of `generic' in each case by using an appropriate `geometric constraint matrix' whose entries are polynomials in the coordinates 
that correspond to the vertices of the underlying combinatorial structure.\\\indent
The results of this paper also immediately apply to realizations of (multi- or hyper-) graphs in spherical $d$-space, since we may interpret a symmetry group of a spherical structure as a symmetry group in Euclidean $(d+1)$-space.\\\indent
Moreover, since infinitesimal and static rigidity are projectively invariant (see \cite{CW, W8, W9}, for example), it is natural to assign projective coordinates instead of Euclidean coordinates to the vertices of the underlying combinatorial structure and to replace a symmetry group consisting of isometries of Euclidean space with a group of projective transformations. We then again obtain linear equations analogous to the ones in Definition \ref{classif}, and hence the definitions and results of this paper (including the definition of `generic') can again immediately be extended to this more general setting.\\\indent Likewise, our results can also be extended to symmetric structures in hyperbolic space.\\\indent
The above-mentioned generalizations allow us to model many different kinds of symmetric structures and to analyze them from a symmetry perspective with regard to questions in various theories related to geometric constraint systems, such as parallel drawings (\cite{W1,W2}), scene analysis (\cite{W1,W2,W7}), or rigidity of body-bar and body-hinge frameworks (\cite{W5, gsw}), for example.\\\indent
This paper therefore provides the foundation for applying symmetry-based techniques to problems not only relating to infinitesimal rigidity of bar and joint frameworks, but also to a range of other geometric constraint problems in various different settings.

\end{document}